\documentclass{amsart}
\usepackage{graphicx,verbatim,amssymb,psfrag}

\newtheorem{proposition}{Proposition}[section]
\newtheorem{theorem}[proposition]{Theorem}

\newtheorem{lemma}[proposition]{Lemma}
\newtheorem{lem}[proposition]{Lemma}
\newtheorem{prop}[proposition]{Proposition}
\newtheorem{cor}[proposition]{Corollary}

\newtheorem{weaker}{Weaker Assertion}

\theoremstyle{definition}
\newtheorem{example}[proposition]{Example}

\theoremstyle{remark}
\newtheorem{remark}[proposition]{Remark}

\numberwithin{equation}{section}

\usepackage[usenames]{color}

\newcounter{margincounter}
\newcommand{\reals}{\mathbb R}

\newcommand{\R}{{\mathcal R}}
\renewcommand{\S}{{\mathcal S}}
\renewcommand{\H}{{\mathcal H}}
\newcommand{\Q}{{\mathcal Q}}
\newcommand{\Po}{{\mathcal P}}

\newcommand{\nc}{\operatorname{nc}}

\newcommand{\Des}{\operatorname{Des}}
\newcommand{\Lower}{\operatorname{Lower}}

\newcommand{\NC}{\operatorname{NC}}
\newcommand{\MC}{\operatorname{MC}}
\newcommand{\rank}{\operatorname{rank}}
\newcommand{\Fix}{\operatorname{Fix}}
\newcommand{\codim}{\operatorname{codim}}

\newcommand{\Int}{\operatorname{Int}}

\newcommand{\Con}{\mbox{{\rm Con}}}
\newcommand{\Can}{\operatorname{Can}}
\newcommand{\Lin}{\operatorname{Lin}}
\newcommand{\Star}{\operatorname{Star}}
\newcommand{\depth}{\operatorname{depth}}

\newcommand{\cov}{\mathrm{cov}}

\newcommand{\covered}{{\,\,<\!\!\!\!\cdot\,\,\,}}
\newcommand{\set}[1]{{\left\lbrace #1 \right\rbrace}}
\newcommand{\pidown}{\pi_\downarrow}
\newcommand{\piup}{\pi^\uparrow}

\newcommand{\br}[1]{{\langle #1 \rangle}}
\newcommand{\A}{{\mathcal A}}
\newcommand{\B}{{\mathcal B}}
\newcommand{\F}{{\mathcal F}}
\newcommand{\K}{{\mathcal K}}
\newcommand{\join}{\vee}
\newcommand{\meet}{\wedge}

\renewcommand{\Join}{\bigvee}
\newcommand{\Meet}{\bigwedge}
\newcommand{\lleq}{\le\!\!\!\le}

\begin{document}
\title[The shard intersection order]{Noncrossing partitions and the shard intersection order}

\author{Nathan Reading}
\address{
Department of Mathematics, North Carolina State University, Raleigh, NC, USA}
\email{nathan\_reading@ncsu.edu}

\begin{abstract}
We define a new lattice structure $(W,\preceq)$ on the elements of a finite Coxeter group~$W.$ 
This lattice, called the \emph{shard intersection order}, is weaker than the weak order and has the noncrossing partition lattice $\NC(W)$ as a sublattice.  
The new construction of $\NC(W)$ yields a new proof that $\NC(W)$ is a lattice. 
The shard intersection order is graded and its rank generating function is the $W$-Eulerian polynomial. 
Many order-theoretic properties of $(W,\preceq)$, like M\"{o}bius number, number of maximal chains, etc., are exactly analogous to the corresponding properties of $\NC(W)$.  
There is a natural dimension-preserving bijection between simplices in the order complex of $(W,\preceq)$ (i.e.\ chains in \mbox{$(W,\preceq)$}) and simplices in a certain pulling triangulation of the $W$-permutohedron.
Restricting the bijection to the order complex of $\NC(W)$ yields a bijection to simplices in a pulling triangulation of the $W$-associahedron.

The lattice $(W,\preceq)$ is defined indirectly via the polyhedral geometry of the reflecting hyperplanes of~$W\!$.  
Indeed, most of the results of the paper are proven in the more general setting of simplicial hyperplane arrangements.
\end{abstract}

\maketitle

\tableofcontents

\section{Introduction}\label{intro}

The (classical) noncrossing partitions were introduced by Kreweras in~\cite{Kreweras}.
Work of Athanasiadis, Bessis, Biane, Brady, Reiner and Watt~\cite{Ath-Rei,Bessis,Biane1,BWKpi,Rei} led to the 
recognition that the classical noncrossing partitions are a special case ($W=S_n$) of a combinatorial construction which yields a noncrossing partition lattice $\NC(W)$ for each finite Coxeter group~$W\!$.

Besides the interesting algebraic combinatorics of the $W$-noncrossing partition lattice, there is a strong motivation for this definition arising from geometric group theory.
In that context, $\NC(W)$ is a tool for studying the Artin group associated to~$W\!$.
(As an example, the Artin group associated to $S_n$ is the braid group.)
For the purposes of Artin groups, a key property of $\NC(W)$ is the fact that it is a lattice.
This was first proved uniformly (i.e.\ without a type-by-type check of the classification of finite Coxeter groups) by Brady and Watt~\cite{BWlattice}.
Another proof, for crystallographic~$W$, was later given by Ingalls and Thomas~\cite{IngThom}.

The motivation for the present work is a new construction of $\NC(W)$ leading to a new proof that $\NC(W)$ is a lattice. 
The usual definition constructs $\NC(W)$ as an interval in a non-lattice (the absolute order) on~$W$; we define a new lattice structure $(W,\preceq)$ on all of~$W$ and identify a sublattice of $(W,\preceq)$ isomorphic to $\NC(W)$.
No part of this construction---other than proving that the sublattice is isomorphic to $\NC(W)$---relies on previously known properties of $\NC(W)$.
Thus, one can take the new construction as a definition of $\NC(W)$.
The proof that $\NC(W)$ can be embedded as a sublattice of $(W,\preceq)$ draws on nontrivial results about \emph{sortable elements} established in \cite{sortable,sort_camb,camb_fan,typefree}.

Beyond the initial motivation for defining $(W,\preceq)$---to construct $\NC(W)$ and prove that it is a lattice---the lattice $(W,\preceq)$ turns out to have very interesting properties.
In particular, many of the properties of $(W,\preceq)$ are precisely analogous to the properties of $\NC(W)$.

The lattice $(W,\preceq)$ is defined in terms of the polyhedral geometry of \emph{shards}, certain \mbox{codimension-1} cones introduced and studied in~\cite{hyperplane,hplanedim,congruence,sort_camb}.
Shards were used to give a geometric description of lattice congruences of the weak order.
In this paper, we consider the collection $\Psi$ of arbitrary intersections of shards, which forms a lattice under reverse containment.
Surprisingly, $\Psi$ is in bijection with~$W\!$.
The lattice $(W,\preceq)$ is defined to be the partial order induced on~$W,$ via this bijection, by the lattice $(\Psi,\supseteq)$.
Thus we call $(W,\preceq)$ the \emph{shard intersection order} on~$W$.

Except in Section~\ref{nc sec}, which deals specifically with $\NC(W)$, most arguments in this paper are given, not in terms of Coxeter groups, but in the slightly more general setting of simplicial hyperplane arrangements.
Although the motivation for this paper lies squarely in the realm of Coxeter groups, it is very natural to argue in the more general setting, because the arguments do not use the group structure of the Coxeter groups at all.
Instead, they rely on the polyhedral geometry of the Coxeter arrangement (a simplicial hyperplane arrangement associated to~$W$) and the lattice structure of weak order on~$W$.

We now summarize the main results in the special case of Coxeter groups.
These are proved later in the paper in the generality of simplicial arrangements, and we indicate, for each result, where to find the more general statement and proof.
Additional results in the body of the paper are phrased only in the broader generality.
In the following propositions, the right descents of $w\in W$ are the simple generators $s\in S$ such that $\ell(ws)<\ell(w)$.
For the proofs, see Propositions~\ref{graded} and~\ref{lower int} and an additional argument given immediately after the proof of Proposition~\ref{lower int}.

\begin{prop}\label{W graded}
The lattice $(W,\preceq)$ is graded, with the rank of $w\in W$ equal to the number of right descents of $w$. 
Alternately, the rank of a cone $C\in\Psi$ is the codimension of~$C$.
\end{prop}
In particular, the rank generating function of $\Psi$ is the $W$-Eulerian polynomial.
For more information on the $W$-Eulerian polynomial, see~\cite{ReinerEul}.

\begin{prop}\label{W lower int}
For any $w\in W$, the lower interval $[1,w]$ in $(W,\preceq)$ is isomorphic to $(W_J,\preceq)$, where $J$ is the set of right descents of $w$ and $W_J$ is the standard parabolic subgroup generated by $J$.
\end{prop}

The identity element of~$W$ is the unique minimal element of $(W,\preceq)$ and the longest element $w_0$ is the unique maximal element.

\begin{theorem}\label{W mobius}
The M\"{o}bius number of $(W,\preceq)$ is $(-1)^{\rank(W)}$ times the number of elements of~$W$ not contained in any proper standard parabolic subgroup.
Equivalently, by inclusion/exclusion,
\[\mu_\preceq(1,w_0)=\sum_{J\subseteq S}(-1)^{|J|}\left|W_J\right|.\]
\end{theorem}
Theorem~\ref{W mobius} is a special case of Theorem~\ref{mobius}.
(An alternate, Coxeter-theoretic proof appears after the proof of Theorem~\ref{mobius}.)
Theorem~\ref{W mobius} is very interesting in light of an analogous description (Theorem~\ref{nc mobius}, due to \cite{ABMW,ABW}) of the M\"{o}bius number of the noncrossing partition lattice $\NC(W)$.
When~$W$ is the symmetric group, the number in Theorem~\ref{W mobius} is, up to sign, the number of indecomposable permutations, or the number of permutations with no global descents.
The latter play a role in the Malvenuto-Reutenauer Hopf algebra of permutations \cite{AgSo}.
See Sequence A003319 in~\cite{Sloane} and the accompanying references.
The corresponding sequences for~$W$ of type $B_n$ or $D_n$ are A109253 and A112225 respectively.

Let $\MC(W)$ be the number of maximal chains in $(W,\preceq)$.
For each $s\in S$, let $\br{s}$ denote $S\setminus\set{s}$.
The following result, proved near the end of Section~\ref{prop sec}, is the only main result on $(W,\preceq)$ without a useful generalization to simplicial arrangements.

\begin{prop}\label{num max}
For any finite Coxeter group~$W$ with simple generators $S$,
\[\MC(W)=\sum_{s\in S}\left(\frac{\left|W\right|}{\left|W_\br{s}\right|}-1\right)\MC(W_\br{s}).\]
\end{prop}

The notation $\MC(W)$ clashes with the author's use (in~\cite{rotate}) of $\MC(W)$ to denote the number of maximal chains in the noncrossing partition lattice $\NC(W)$.
In fact, the number of maximal chains of $\NC(W)$ satisfies a recursion \cite[Corollary~3.1]{rotate} very similar to Proposition~\ref{num max}.
The latter recursion can be solved non-uniformly \cite[Theorem~3.6]{rotate} to give a uniform formula first pointed out in \cite[Proposition~9]{Chapoton}.

Recursions involving sums over maximal proper parabolic subgroups, such as the recursion appearing in Proposition~\ref{num max}, are very natural in the context of Coxeter groups and root systems.
Besides Proposition~\ref{num max} and \cite[Corollary~3.1]{rotate}, there are at least two other important examples:
One is a a recursive formula for the face numbers of generalized associahedra \cite[Proposition~3.7]{ga}. 
(Cf. \cite[Proposition~8.3]{gcccc}.)
Yet another is a formula for the volume of the~$W$-permutohedron which can be obtained by simple manipulations from Postnikov's formula \cite[Theorem~18.3]{Postnikov} expressing volume in terms of $\Phi$-trees.

Let $\Delta(W)$ be the pulling triangulation of the $W$-permutohedron, where the vertices are ordered by the reverse of the weak order.
This construction is described in more detail in Section~\ref{prop sec};  see also~\cite{Lee}.
\begin{theorem}\label{perm chain tri}
There exists a dimension-preserving bijection between $\Delta(W)$ and the order complex of $(W,\preceq)$.
\end{theorem}
In particular, the $f$-vector of the order complex of $(W,\preceq)$ coincides with the $f$-vector of $\Delta(W)$.
For the proof, see Theorem~\ref{zon chain tri}.

Since $(W,\preceq)$ is defined in terms of shards, which encode lattice congruences of the weak order, it should not be surprising that $(W,\preceq)$ is compatible with lattice congruences on the weak order.
Specifically, given a lattice congruence~$\Theta$ on the weak order, let $\pidown^\Theta(W)$ denote the set of minimal-length congruence class representatives.
The restriction $(\pidown^\Theta(W),\preceq)$ of the shard intersection order to $\pidown^\Theta(W)$ is a join-sublattice of $(W,\Theta)$ and shares many of the properties of $(W,\preceq)$.
In particular, direct generalizations of Propositions~\ref{W graded} and~\ref{W lower int} and Theorems~\ref{W mobius} and~\ref{perm chain tri} are stated and proved in Section~\ref{cong psi sec}.

For each Coxeter element~$c$ of~$W$, there is a noncrossing partition lattice $\NC_c(W)$.
The isomorphism type of $\NC_c(W)$ is independent of $c$, so we suppressed the dependence on~$c$ earlier in the introduction.
The \emph{$c$-Cambrian congruence}~$\Theta_c$ is a lattice congruence defined in~\cite{cambrian} and studied further in \cite{sort_camb,camb_fan,typefree}.
The set $\pidown^{\Theta_c}(W)$  can be characterized combinatorially~\cite{sort_camb} as the set of \emph{$c$-sortable elements} of~$W$.
As a special case of general results from~\cite{con_app}, the $c$-Cambrian lattice defines a complete fan which coarsens the fan defined by the Coxeter arrangement $\A(W)$ of reflecting hyperplanes of~$W$.
This fan is combinatorially isomorphic~\cite{camb_fan} to the normal fan of the $W$-associahedron, which was defined in~\cite{ga}.

Our discussion of noncrossing partitions is found in Section~\ref{nc sec}.
As a special case of the general result mentioned above, the $c$-sortable elements induce a join-sublattice $(\pidown^{\Theta_c}(W),\preceq)$ of $(W,\preceq)$.
Drawing on results of \cite{sortable,typefree}, we show that $(\pidown^{\Theta_c}(W),\preceq)$ is isomorphic to $\NC_c(W)$.
In particular, we obtain not only a new proof of the lattice property for $\NC_c(W)$ but also a completely new construction of $\NC_c(W)$.
Furthermore, we show that $(\pidown^{\Theta_c}(W),\preceq)$ is a sublattice of $(W,\preceq)$, rather than merely a join-sublattice.
Applying general results on $(\pidown^{\Theta}(W),\preceq)$ to the case $\Theta=\Theta_c$, we give new proofs of old and new results on noncrossing partitions.
In particular, we generalize, to arbitrary~$W$, a bijection of Loday~\cite{Loday} from maximal chains in the classical noncrossing partition lattice (in the guise of parking functions) to maximal simplices in a certain pulling triangulation of the associahedron.
We also broaden the bijection into a dimension-preserving bijection (Theorem~\ref{chain tri camb}) between simplices in the order complex of $\NC_c(W)$ and simplices in the triangulation.
Both Theorem~\ref{chain tri camb} and the analogous statement (Theorem~\ref{perm chain tri}) for the permutohedron and $(W,\preceq)$  are special cases of a much more general result, Theorem~\ref{chain tri cong}.

The construction of noncrossing partitions via shard intersections exhibits a surprising connection to semi-invariants of quivers, which we hope to explain more fully in a future paper.
Some additional detail is given in Remark~\ref{semiinvariant remark}.

\section{Simplicial hyperplane arrangements}\label{arr sec}
This section covers background information on simplicial hyperplane arrangements that is used in the rest of the paper.
We also explain how the weak order on a finite Coxeter group fits into the context of simplicial hyperplane arrangements.

A \emph{linear hyperplane} in a vector space~$V$ is a codimension-1 linear subspace of~$V$.
An \emph{affine hyperplane} in~$V$ is any translate of a linear hyperplane.
A \emph{hyperplane arrangement}~$\A$ in~$V$ is a finite collection of hyperplanes.
Without exception, throughout the paper, we take~$\A$ to be \emph{central}, meaning that all hyperplanes in~$\A$ are linear.

The \emph{rank} of a central arrangement~$\A$ is the codimension of the intersection $\bigcap\A$ of all the hyperplanes in~$\A$.
A central hyperplane arrangement~$\A$ is called \emph{essential} if $\bigcap\A$ has dimension zero.
We do not require our arrangements to be essential, because it is convenient to consider an arrangement~$\A$ in the same vector space as a subarrangement $\A'\subseteq\A$, even when $\A'$ has lower rank.
However, it is easy to make an essential arrangement with the same combinatorial structure as~$\A$ by passing to the quotient vector space $V/(\bigcap\A)$, and thus the reader may safely think in the essential case.
A central hyperplane arrangement~$\A$ is a \emph{direct sum} of $\A_1$ and $\A_2$ if $\A=\A_1\cup\A_2$ and $V$ is a direct sum $V=V_1\oplus V_2$ such that $\A_1=\set{H\in\A:H^\perp\in V_1}$ and $\A_2=\set{H\in\A:H^\perp\in V_2}$.

A \emph{region} of~$\A$ (or ``$\A$-region'') is the closure of a connected component of the complement $V\setminus(\bigcup\A)$ of~$\A$.
Each region of a central arrangement is a closed convex polyhedral cone whose dimension equals $\dim(V)$.
(A convex polyhedral cone is a set of points determined by a finite system of linear inequalities.)
The set of regions is denoted by~$\R$ or $\R(\A)$.
We speak of \emph{faces} of a region~$R$ in the usual polyhedral sense.
A \emph{facet} of~$R$ is a maximal proper face of~$R$.
A region is \emph{simplicial} if the normal vectors to its facet-defining hyperplanes form a linearly independent set.
When~$\A$ is essential, a region is simplicial if and only if it is a cone over a simplex.
A central hyperplane arrangement~$\A$ is \emph{simplicial} if every $\A$-region is simplicial.

We now fix a \emph{base region} $B\in\R$, and define the \emph{poset of regions} $(\R(\A),\le_B)$ or simply $(\R,\le)$.
(In \cite{BEZ,hyperplane,hplanedim,congruence,con_app}, this poset is denoted by $\Po(\A,B)$ or $\Po(\H,B)$.)
Given a region $R\in\R$, the \emph{separating set} $S(R)$ of~$R$ is the set of hyperplanes $H\in\A$ such that~$H$ separates~$R$ from~$B$.
The poset of regions sets $Q\le R$ if and only if $S(Q)\subseteq S(R)$.
This partially ordered set is a lattice \cite[Theorem~3.4]{BEZ} with a unique minimal element~$B$ and a unique maximal element $-B$.
Cover relations in $(\R,\le)$ are $Q\covered R$ such that $Q$ and~$R$ share a facet, and the hyperplane defining that facet separates~$R$ from~$B$.
The involution $R\mapsto(-R)$ is an anti-automorphism of $(\R,\le)$.

Let $(W,S)$ be a finite Coxeter group, and represent~$W$ in the usual way as a group of orthogonal transformations of some Euclidean vector space~$V$.
The set $T=\set{wsw^{-1}:w\in W,s\in S}$ is the collection of all elements of~$W$ that act as reflections in~$V$.
For each reflection $t\in T$, let~$H_t$ be the hyperplane fixed by~$t$.
The \emph{Coxeter arrangement} $\A(W)=\set{H_t:t\in T}$ is a central, simplicial hyperplane arrangement whose rank equals the rank of~$W$.
The regions of $\A(W)$ are in bijection with the elements of~$W$ as follows:  one first chooses the base region~$B$ to be a region whose facet-defining hyperplanes are $\set{H_s:s\in S}$.
(There are two choices, related by the antipodal map.)
The bijection maps $w\in W$ to the region $wB$.

\begin{figure}
\begin{tabular}{ccccc}
\scalebox{.9}{\includegraphics{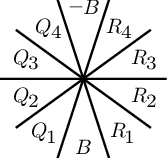}}&&
\scalebox{.9}{\includegraphics{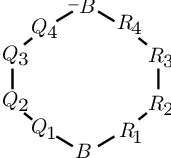}}&&
\scalebox{.9}{\includegraphics{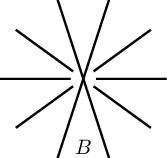}}\\[2 pt]
(a)&&(b)&&(c)
\end{tabular}
\caption{a: An arrangement of five lines in $\reals^2$.  b: The corresponding poset of regions. c: The shards.}
\label{i25 fig}
\end{figure}

\begin{example}\label{labeled i25}
Let~$\A$ be a set of $5$ distinct lines through the origin in $\reals^2$.
If the lines all meet at equal angles, then~$\A$ is the Coxeter arrangement for the dihedral Coxeter group $I_2(5)$.
Figure~\ref{i25 fig}.a shows the arrangement~$\A$ with the 10 regions labeled.
Figure~\ref{i25 fig}.b shows the poset of regions $(\R(\A),\le_B)$.
\end{example}

\begin{figure}
\scalebox{.85}{\includegraphics{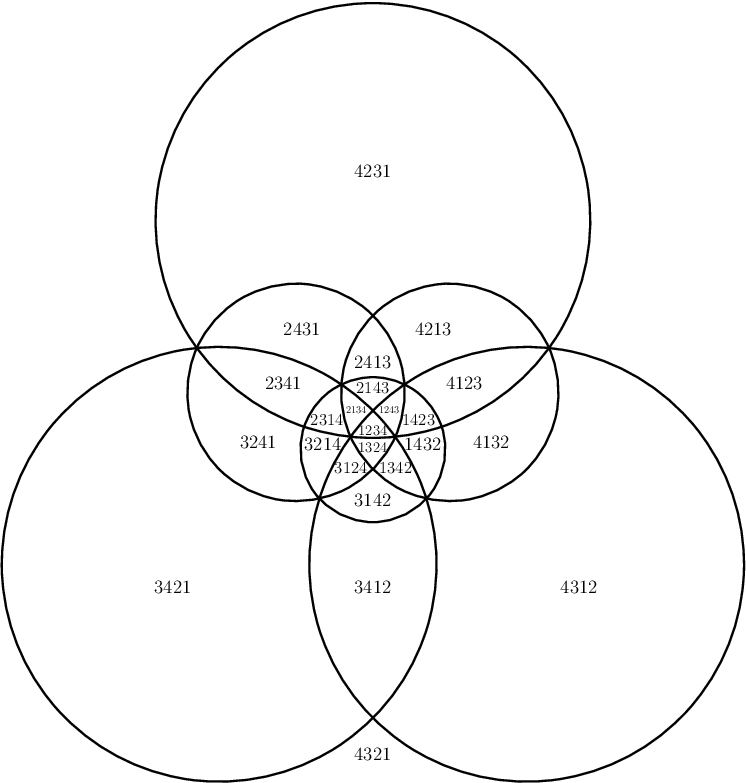}}
\caption{The Coxeter arrangement $\A(W)$ for $W=S_4$.}
\label{A3labels}
\end{figure}

\begin{example}\label{labeled S4}
The Coxeter group~$W$ of type $A_3$ is isomorphic to the symmetric group~$S_4$.
The Coxeter arrangement~$\A(W)$ consists of six hyperplanes through the origin in $\reals^3$.
These planes, intersected with the unit sphere in $\reals^3$, define an arrangement of six great circles on the sphere.
A stereographic projection yields an arrangement of six circles in the plane.
This arrangement of circles is shown in Figure~\ref{A3labels}.
Regions of~$\A$ appear as curved-sided triangles.  
Each region is labeled with the corresponding permutation in $S_4$.
We choose the base region~$B$ to be the small triangle, labeled $1234$, which is inside the three large circles.
Of necessity, some labels near the center of the picture are quite small.
These are included for the benefit of readers viewing this paper electronically.
For the benefit of those reading this paper in print:
The label on the triangle inside all circles is $1234$.
The label on its lower neighbor is $1324$, the label on its top-left neighbor is $2134$ and the label on its top-right neighbor is $1243$. 
\end{example}

The weak order is a partial order on~$W$ which can be defined combinatorially in terms of reduced words.
(There are two isomorphic weak orders on $W$; we consider the ``right'' weak order, as opposed to the ``left'' weak order.)
Alternately, the weak order is defined in terms of containment of inversion sets.
In the latter guise, the weak order is seen to be isomorphic to the poset of regions $(\R(\A(W)),\le_B)$, for the choice of~$B$ described above.

\begin{figure}
\scalebox{.75}{\includegraphics{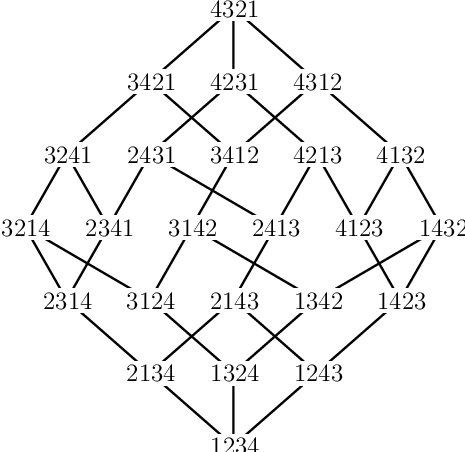}}
\caption{The weak order on $S_4$.}
\label{A3weak}
\end{figure}

\begin{example}\label{weak ex}
The weak order on $S_n$ can be described in terms of the combinatorics of permutations.
A covering pair consists of two permutations which agree, except that two adjacent entries have been swapped.
The lower permutation of the two, in the weak order, is the permutation in which the two adjacent entries occur in numerical order.
The weak order on $S_4$ is shown in Figure~\ref{A3weak}.
\end{example}

For the remainder of the paper, we assume that $\A$ is a simplicial hyperplane arrangement with a chosen base region~$B$.
Let $\B$ be the set of facet-defining hyperplanes of~$B$.
Since~$\A$ is simplicial, the cardinality of $\B$ is equal to the rank of~$\A$.
Given a set $\K\subseteq \B$, let $\A_\K$ be the set $\set{H\in\A:H\supseteq (\bigcap \K)}$.
The arrangement $\A_\K$ is called a \emph{standard subarrangement} of~$\A$.
(In \cite{congruence}, the term \emph{parabolic subarrangement} was used for what we are here calling a standard subarrangement of~$\A$.)
Also associated to $\K$ is a subset $\R_\K$ of~$\R$ defined as follows:
For each $\A$-region~$R$, there exists \cite[Lemma~6.2]{congruence} a (necessarily unique) $\A$-region $R_\K$ such that $S(R_\K)=S(R)\cap\A_\K$.
The set $\R_\K=\set{R_\K:R\in\R}$, called a \emph{standard parabolic subset} of~$\R$, is the set of regions whose separating sets are contained in the standard subarrangement $\A_\K$.

Standard subarrangements are a special case of a more general notion.
If $\A'\subseteq\A$ is the collection of all hyperplanes in~$\A$ containing a particular subset of~$V$ then $\A'$ is called a \emph{full subarrangement} of~$\A$.
Let $B'$ be the $\A'$-region containing~$B$.
The \emph{basic hyperplanes} of $\A'$ are the facet-defining hyperplanes of $B'$.
The set of basic hyperplanes of~$\A$ is $\B$ and the set of basic hyperplanes of a standard subarrangement $\A_\K$ is $\K$.

Full subarrangements of a Coxeter arrangement correspond to parabolic subgroups of the Coxeter group.
The subarrangement is the set of reflecting hyperplanes of the parabolic subgroup;  the parabolic subgroup is the subgroup generated by reflections in hyperplanes of the subarrangement.
In the same sense, standard subarrangements of a Coxeter arrangement correspond to standard parabolic subgroups.
Standard parabolic subsets also correspond to standard parabolic subgroups of~$W$, but in a different sense:
Recall that the set~$\R$ of $\A(W)$-regions is in bijection with the elements of $W$; a standard parabolic subset $\R_\K$ is the set of $\A(W)$-regions corresponding to elements of $W_K$, where $K=\set{s\in S:H_s\in\K}$.

\begin{example}\label{subarrangement}
This example refers to Figure~\ref{A3labels}, which represents the Coxeter arrangement $\A(W)$ for $W=S_4$, as explained in Example~\ref{labeled S4}.
The basic hyperplanes $\B$ of $\A(W)$ are represented by the circles defining the boundary of the regions~$B$ (labeled by $1234$).
These are~$H_{(1\,2)}$,~$H_{(2\,3)}$ and~$H_{(3\,4)}$, the hyperplanes separating~$B$ from the regions labeled $2134$, $1324$ and $1243$.
Consider the point~$p$ defined as the intersection of the triangle labeled $3421$ with the triangle labeled $3124$.
The set of three circles containing~$p$ describes a (nonstandard) full subarrangement $\A'$ of $\A(S_4)$.
(There is a ray in $\reals^3$ whose projection is~$p$, and $\A'$ is the set of hyperplanes containing that ray.)
The basic hyperplanes of $\A'$ are represented by the circle separating $3124$ from $3214$ and the circle separating $3124$ from $3142$.
Consider $\K=\set{H_{(1\,2)},H_{(2\,3)}}\subseteq\B$.
The standard subarrangement $\A(W)_K$ consists of all of the hyperplanes of $\A(W)$ containing the intersection of~$H_{(1\,2)}$ and~$H_{(2\,3)}$.
These are the three hyperplanes separating the region labeled $3214$ from the region~$B$, labeled $1234$.
The standard parabolic subset $\R_\K$ of~$\R$ consists of regions labeled $\set{1234,2134,1324,2314,3124,3214}$.
\end{example}

The following lemma will be useful in the next section.
\begin{lemma}\label{parabolic lem}
Let $\K\subseteq \B$, let~$H_1\in(\A-\A_\K)$, let~$H_2\in\A_\K$ and let $\A'$ be the rank-two full subarrangement containing~$H_1$ and~$H_2$.
Then $(\A'\cap\A_\K)=\set{H_2}$ and~$H_2$ is basic in $\A'$.
\end{lemma}
\begin{proof}
The special case where $|\K|=|\B|-1$ is precisely the statement of \cite[Lemma~6.6]{congruence}.
The general result follows easily from the fact that any subset $\K\subseteq\B$ is the intersection of subsets $\K'\subseteq \B$ with $|\K'|=|\B|-1$, together with the observation that $\A_{\K'\cap\K''}=\A_{\K'}\cap\A_{\K''}$ for any $\K',\K''\subseteq\B$.
\end{proof}

The collection of regions, together with all of their faces, forms a complete fan~$\F$ or $\F(\A)$ in~$V$.
(For background information on fans, see, for example, \cite[Lecture~7]{Ziegler}.)
The interaction of the fan~$\F$ with the poset of regions $(\R,\le)$ is discussed extensively in~\cite{con_app}.
The following proposition and theorem summarize a very small part of the discussion, found at the beginning of \cite[Section~4]{con_app}.

\begin{prop}\label{facial}
For any face $F\in\F$, the set $\set{P\in \R:P\supseteq F}$ is an interval $[Q,R]$ in $(\R,\le)$.
Furthermore,~$F$ is the intersection of the facets of~$R$ separating~$R$ from cones $P$ with $Q\le P\covered R$.
Dually,~$F$ is the intersection of the facets of $Q$ separating $Q$ from cones $P$ with $Q\covered P\le R$.
The interval $[Q,P]$ is isomorphic to the poset of regions $(\R(\A'),\le_{B'})$, where $\A'$ is the full subarrangement of~$\A$ consisting of hyperplanes in~$\A$ containing~$F$, and $B'$ is the $\A'$-region containing~$B$.
\end{prop}

\begin{theorem}\label{shelling}
Any linear extension of the poset of regions $(\R,\le)$ (or of its dual $(\R,\ge)$) is a shelling order on the maximal cones of $\F(\A)$.
\end{theorem}
In Theorem~\ref{shelling}, the assertion about $(\R,\le)$ is, by \cite[Proposition~4.2]{con_app}, a special case of a more general result, \cite[Proposition~3.4]{con_app}.
The assertion about $(\R,\ge)$ follows because $R\mapsto-R$ is an antiautomorphism of $(\R,\le)$ and the antipodal map is an automorphism of~$\F$.
In the following lemma, $\bigcap\emptyset=V$ by convention.

\begin{lemma}\label{RK}
If $\K$ is a subset of the basic hyperplanes of~$\A$ then $\R_\K$ is the set of regions containing the face $B\cap\bigcap\K$ of~$\F$.
Furthermore, $|\R_\K|$ coincides with the number of regions containing the face $(-B)\cap\bigcap\K$ of~$\F$.
\end{lemma}
\begin{proof}
Suppose~$R$ contains the face $F=B\cap\bigcap\K$.
Then one can choose a point $x$ in the interior of~$B$ and $y$ in the interior of~$R$ such that the line segment $\overline{xy}$ intersects~$F$ and no other face of~$B$ or of~$R$.
Thus $S(R)$ contains only hyperplanes containing~$F$, or equivalently, only hyperplanes containing $\bigcap\K$.
In other words, $R=R_\K$.
This argument is easily reversed.

The second assertion follows because $R\mapsto-R$ is an involution on~$\R$ and the antipodal map is an automorphism of~$\F$.
\end{proof}

A point $x$ in~$V$ is said to be \emph{below} a hyperplane $H\in\A$ if $x$ is contained in~$H$ or if $x$ and~$B$ are on the same side of~$H$.
The point $x$ is \emph{strictly below}~$H$ if $x$ is below~$H$ but is not contained in~$H$.
A subset of~$V$ is below, or strictly below~$H$ if each of its points is.
The notions of \emph{above} and \emph{strictly above} are defined similarly.

In the same spirit, given a region $R\in\R$, we define a \emph{lower hyperplane} of~$R$ to be a hyperplane in~$\A$ containing a facet of~$R$ which separates~$R$ from a region $Q\covered R$.
The set of lower hyperplanes of~$R$ is written $\Lower(R)$.

\begin{prop}\label{cov ref basic}
For any $R\in\R$, let~$F$ be the intersection of all facets of~$R$ separating~$R$ from a region covered by~$R$ and let $\A'$ be the full subarrangement consisting of hyperplanes containing~$F$.
Then the lower hyperplanes $\Lower(R)$ of~$R$ are the basic hyperplanes of $\A'$.
\end{prop}
\begin{proof}
Let $B'$ be the $\A'$-region containing~$B$.
Then the basic hyperplanes of $\A'$ are the facet-defining hyperplanes of $B'$, or equivalently, the facet-defining hyperplanes of $-B'$.
These coincide with the facet-defining hyperplanes of~$R$ that contain~$F$, or equivalently, the lower hyperplanes $\Lower(R)$ of~$R$.
\end{proof}

We conclude the section with a useful technical lemma.

\begin{lemma}\label{tech}
Let~$\F$ be a complete fan of convex polyhedral cones.
Let $C_1$ and $C_2$ be convex polyhedral cones, each of which is a union of faces of~$\F$.
Let $F_1$ be a face of~$\F$ contained in $C_1$ with $\dim(F_1)=\dim(C_1)$.
If $C_1\subseteq C_2$ then there exists a face $F_2$ of~$\F$ with $F_1\subseteq F_2\subseteq C_2$ and $\dim(F_2)=\dim(C_2)$.
\end{lemma}
\begin{proof}
Let $x$ be a vector in the relative interior of $F_1$ and let $y$ be a generic vector in the relative interior of $C_2$.
For sufficiently small positive $\epsilon$, the vector $(1-\epsilon)x+\epsilon y$ is in the relative interior of a face $F_2$ with the desired properties.
\end{proof}

\section{Cutting hyperplanes into shards}\label{shard sec}
Recall that, throughout the paper,~$\A$ is a simplicial hyperplane arrangement and~$B$ is a choice of base region. 
In the first half of this section, we review a ``cutting'' relation on hyperplanes in~$\A$ and review the use of the cutting relation to define shards.
Proposition~\ref{shard ji} and Theorem~\ref{R can join rep}, below, provide some motivation for the notion of shards.
Further motivation for shards, arising from the study of lattice congruences on $(\R,\le)$, appears later in Section~\ref{cong sec}.
The second part of this section is devoted to proving lemmas which are crucial in the study of intersections of shards.
Although the definition of shards is valid even when~$\A$ is not simplicial, most results discussed in this section rely on the assumption that~$\A$ is simplicial.

The cutting relation depends implicitly on the choice of~$B$.
Given two hyperplanes $H,H'$ in~$\A$, let $\A'(H,H')$ be the full subarrangement of~$\A$ consisting of hyperplanes of~$\A$ containing $H\cap H'$.
The subarrangement $\A(H,H')$ has rank two, and is in fact the unique rank-two full subarrangement containing $H$ and $H'$.
We say that~$H$ \emph{cuts} $H'$ if~$H$ is a basic hyperplane of $\A(H,H')$ and $H'$ is \textbf{not} a basic hyperplane of $\A(H,H')$.
For each $H\in\A$, remove from~$H$ all points contained in hyperplanes of~$\A$ that cut~$H$.
The remaining set of points may be disconnected; the closures of the connected components are called the \emph{shards} in~$H$.
Thus~$H$ is ``cut'' into shards by certain hyperplanes in~$\A$, just as~$V$ is ``cut'' into regions by all of the hyperplanes in~$\A$.
The set of shards of $\A$ is the union, over hyperplanes $H\in\A$, of the set of shards in $H$.
(In~\cite{hyperplane} and~\cite{hplanedim}, shards were defined to be the relatively open connected components, without taking closures.  All results on shards cited from these sources have been rephrased as necessary.)

\begin{example}\label{i25shards}
This is a continuation of Example~\ref{labeled i25}.
The $8$ shards in the arrangement of 5 lines in $\reals^2$ are illustrated in Figure~\ref{i25 fig}.c.
Each is a one-dimensional cone.
The two lines intersecting at the origin are two distinct shards.
All of the shards contain the origin; however, some shards in the picture are offset slightly to indicate that they do not continue through the origin.
\end{example}

\begin{figure}
\scalebox{1}{\includegraphics{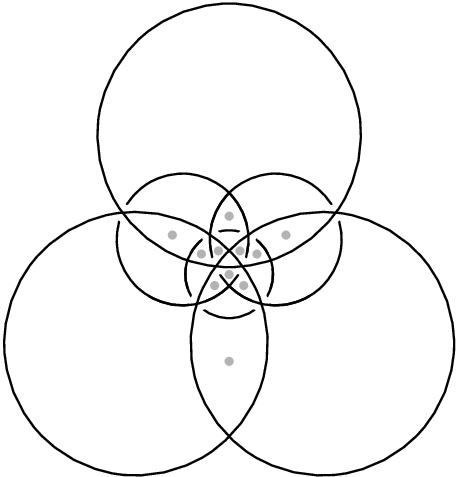}}
\caption{Shards in the Coxeter arrangement $\A(S_4)$.}
\label{A3shards}
\end{figure}

\begin{example}\label{A3shards ex}\rm
The shards in the Coxeter arrangement $\A(W)$, for the case $W=S_4$, are pictured in Figure~\ref{A3shards}.
This figure is a stereographic projection as explained in Example~\ref{labeled S4}.
As before, the cone~$B$ is the small triangular region which is inside the three largest circles.
The shards are closed two-dimensional cones (which in some cases are entire planes).
Thus they appear as full circles or as circular arcs in the figure.
To clarify the picture, we continue the convention of Figure~\ref{i25 fig}.c: 
Where shards intersect, certain shards are offset slightly from the intersection to indicate that they do not continue through the intersection.
Some of the regions are marked with gray dots.
The significance of these regions is explained in Example~\ref{A3 shard ji}.
\end{example}

The unique hyperplane containing a shard~$\Sigma$ is denoted by $H(\Sigma)$.
An \emph{upper region} of a shard~$\Sigma$ is a region $R\in\R$ such that $\dim(R\cap\Sigma)=\dim(\Sigma)$ and $H(\Sigma)\in S(R)$.
That is, a region of~$\A$ is an upper region of~$\Sigma$ if it has a facet contained in~$\Sigma$ such that the region adjacent through that facet is lower (necessarily by a cover) in the poset of regions.
Let $U(\Sigma)$ be the set of upper regions of~$\Sigma$, partially ordered as an induced subposet of the poset of regions.

An element $j$ in a finite lattice~$L$ is called \emph{join-irreducible} if it covers exactly one element, denoted $j_*$.
The following proposition is a concatenation of \cite[Proposition~2.2]{hplanedim} and \cite[Proposition~3.5]{congruence}.

\begin{prop}\label{shard ji}
For any shard~$\Sigma$, there is a unique minimal element of $U(\Sigma)$.
This region, denoted by $J(\Sigma)$, is join-irreducible in $(\R,\le)$, and furthermore every join-irreducible element of $(\R,\le)$ is $J(\Sigma)$ for a unique shard~$\Sigma$.
\end{prop}

\begin{example}\label{A3 shard ji}
The regions corresponding to join-irreducible elements of the weak order on $S_4$ (the poset of regions of $\A(S_4)$) are marked in Figure~\ref{A3shards} by gray dots.
Each dotted triangle has two convex sides and one concave side.
The bijection between join-irreducible regions and shards sends the triangle to the shard containing its concave side.
\end{example}

The notation $\Sigma(J)$ denotes the unique shard~$\Sigma$ such that $J=J(\Sigma)$.
We now give a stronger characterization of $J(\Sigma)$.
If~$R$ is an upper region of~$\Sigma$, then we say that~$\Sigma$ is a \emph{lower shard} of~$R$.

\begin{lemma}\label{J Sig char}
Let~$\Sigma$ be a lower shard of $R\in\R$.
Then $J(\Sigma)$ is the unique minimal region in $(\R,\le)$ among regions $Q\le R$ with $H(\Sigma)\in S(Q)$.
\end{lemma}
\begin{figure}
\scalebox{.65}{\includegraphics{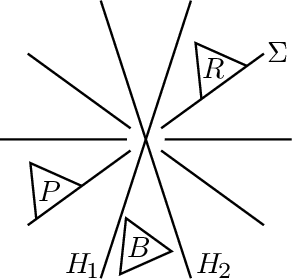}}
\caption{An illustration of the proof of Lemma~\ref{J Sig char}.}
\label{schematic}
\end{figure}
\begin{proof}
Let $J=J(\Sigma)$.
Since $J$ is an upper region of~$\Sigma$, in particular $H(\Sigma)\in S(J)$.
Since~$R$ is also an upper region of~$\Sigma$, Proposition~\ref{shard ji} says that $J\le R$.
If $Q$ is any region with $Q\le R$ and $H(\Sigma)\in S(Q)$, then there exists $P\le Q$ such that $H(\Sigma)$ is a lower hyperplane of $P$.
(To find such a $P$, consider an unrefinable chain in $(\R,\le)$ from~$B$ to $Q$.
Since $H(\Sigma)\in S(Q)$, there exists a covering pair $P'\covered P$ in the chain such that $H(\Sigma)\in S(P)$ but $H(\Sigma)\not\in S(P')$.)

We claim that $P$ is an upper region of~$\Sigma$.
If not, then $P\cap H(\Sigma)$ is separated from~$\Sigma$ by the intersection of $H(\Sigma)$ with a hyperplane that cuts $H(\Sigma)$.
In fact, there are two such hyperplanes,~$H_1$ and~$H_2$ which cut $H(\Sigma)$ in the same place.
Simple geometric considerations (illustrated schematically in Figure~\ref{schematic}) show that, without loss of generality,~$H_1\in S(P)$ and~$H_1\not\in S(R)$.
This contradicts the fact that $P\le Q\le R$, thus proving the claim.
Since $P$ is an upper region of~$\Sigma$, we have $Q\ge P\ge J(\Sigma)$ by Proposition~\ref{shard ji}.
\end{proof}

Any cover relation in $(\R,\le)$ uniquely determines a shard:
Given $Q\covered R$ in $(\R,\le)$, the intersection $Q\cap R$ is a facet of $Q$ and of~$R$.
There is a unique shard containing this facet, denoted by $\Sigma(Q\covered R)$.

We now define the \emph{canonical join-representation} of an element of a finite lattice~$L$.
The canonical join representation of $x\in L$, when it exists, is the set $\Can(x)$ such that $\Join\Can(x)$ is the unique ``lowest'' non-redundant expression for $x$ as a join, in a sense which we now make precise.
An expression $x=\Join A$ is \emph{redundant} if some proper subset $A'\subsetneq A$ has $x=\Join A'$.
The requirement that $\Join\Can(x)$ be a non-redundant expression implies in particular that $\Can(x)$ is an antichain (a set of pairwise incomparable elements) in~$L$.
To define $\Can(x)$, we define a partial order $\lleq$ on antichains in~$L$ by setting $A\lleq B$ if and only if for every $a\in A$ there exists $b\in B$ with $a\le b$.
Then $\Can(x)$ is the unique minimal antichain, with respect to $\lleq$, among antichains joining to $x$, if this unique minimal antichain exists.
The elements of $\Can(x)$ are called \emph{canonical joinands} of $x$.
It is easily checked that every canonical joinand of $x$ is a join-irreducible element of~$L$.
It is also easily seen than a proper subset $A\subsetneq\Can(x)$ is the canonical join-representation of some element $x'<x$.
For more information on canonical join-representations, see \cite[Section~II.1]{FreeLattices}.

The following theorem is essentially \cite[Theorem~8.1]{typefree}.
However, the latter result is more special than Theorem~\ref{R can join rep} because it is proven for the weak order on a Coxeter group, but more general than Theorem~\ref{R can join rep} in that it allows the Coxeter group to be infinite.

\begin{theorem}\label{R can join rep}
Every $R\in\R$ has a canonical join representation in $(\R,\le)$, namely the set of regions $J(\Sigma)$, where~$\Sigma$ ranges over all lower shards of~$R$.
Furthermore $\Lower(R)$ is the disjoint union, over canonical joinands $J$, of the singletons $\Lower(J)$.
\end{theorem}

\begin{proof}
Let the lower shards of~$R$ be $\Sigma_1,\ldots,\Sigma_k$.
Lemma~\ref{J Sig char} (or Proposition~\ref{shard ji}) implies that $R\ge J(\Sigma_i)$ for $i\in[k]$.
On the other hand, any element $Q\covered R$ is separated from~$R$ by a hyperplane $H(\Sigma_i)$, so $H(\Sigma_i)\not\in S(Q)$ and thus $Q\not\ge J(\Sigma_i)$.
Since $(\R,\le)$ is a lattice,~$R$ must be $J(\Sigma_1)\join\cdots\join J(\Sigma_k)$.

For any $i\in[k]$, there is a region $Q_i$ covered by~$R$ which is separated from~$R$ by $H(\Sigma_i)$ and no other hyperplane.
For all $j\in[k]$ with $j\neq i$, we have $H(\Sigma_j)\in S(Q_i)$ and $Q_i\le R$, so Lemma~\ref{J Sig char} implies that $Q_i\ge J(\Sigma_j)$.
We conclude that the join of any proper subset of $\set{J(\Sigma_1),\ldots,J(\Sigma_k)}$ is strictly smaller than~$R$.
Thus $R=\Join\set{J(\Sigma_1),\ldots,J(\Sigma_k)}$ is a non-redundant expression for~$R$ and in particular $\set{J(\Sigma_1),\ldots,J(\Sigma_k)}$ is an antichain in $(\R,\le)$.

Let $A$ be any other antichain in $(\R,\le)$ having $\Join A=R$.
Let $i\in[k]$.
Some element $P_i$ of $A$ has $H(\Sigma_i)\in S(P)$:
Otherwise, the region $Q_i$, defined in the previous paragraph, is an upper bound for $A$, contradicting $\Join A=R$.
Thus $P_i\ge J(\Sigma_i)$ by Lemma~\ref{J Sig char}. 
Now $\set{J(\Sigma_1),\ldots,J(\Sigma_k)}\lleq A$, and we have proved that $\set{J(\Sigma_1),\ldots,J(\Sigma_k)}$ equals $\Can(R)$.
\end{proof}

\begin{example}\label{S4 can join rep}
We give an example of Theorem~\ref{R can join rep}, for the case $\A=\A(S_4)$.
Consider the element $4312\in S_4$.
It is easily verified using Figure~\ref{A3weak} that $4312=3124\join1243$, and that the set $\set{3124,1243}$ is minimal in the order $\lleq$ among antichains joining to $4312$.
Additional inspection of Figure~\ref{A3weak} shows that $\set{3124,1243}$ is the unique $\lleq$-minimal antichain joining to $4312$, or in other words, that $\set{3124,1243}$ is the canonical join representation of $4312$.

Referring to Figure~\ref{A3labels} for the labeling, we can find the lower shards of the region labeled $4312$ in Figure~\ref{A3shards}.
They are the shards containing the concave side of the triangle.
The minimal upper regions of these two shards are the join-irreducible regions which (again referring to Figure~\ref{A3labels}) are labeled $3124$ and $1243$.
\end{example}

The next three lemmas detail the interaction between the cutting relation and full subarrangements.
The first is immediate from the definition, the second follows immediately from Lemma~\ref{parabolic lem} and the definition, and the third is \cite[Lemma 6.8]{congruence}.

\begin{lemma}\label{full cutting}
Let $\A'$ be a full subarrangement of~$\A$, let $B'$ be the $\A'$-region containing~$B$.
Then the cutting relation on $\A'$, defined with respect to $B'$, is the restriction of the cutting relation on~$\A$ to hyperplanes in $\A'$.
\end{lemma}

\begin{lemma}\label{para cutting}
If $\K\subseteq \B$ and $H\in\A_\K$ then~$H$ is not cut by any hyperplane of $\A\setminus\A_\K$.
In particular, if~$\Sigma$ is a shard of~$\A$ contained in a hyperplane $H\in\A_\K$, then $\Sigma\supseteq(\bigcap \K)$.
\end{lemma}

\begin{lemma}\label{para shard}
If $\K\subseteq \B$ and~$\Sigma$ is a shard then $H(\Sigma)\in\A_\K$ if and only if \mbox{$J(\Sigma)\in\R_\K$}.
\end{lemma}

The following lemma is a slight rephrasing\footnote{The equivalence of condition (i) of Lemma~\ref{whole} and condition (i) of \cite[Lemma~3.9]{congruence} is explained in the first paragraph of the proof of \cite[Lemma~3.9]{congruence}.}
 of \cite[Lemma~3.9]{congruence}.
\begin{lemma}
\label{whole}
Let~$\Sigma$ be a shard.
Then the following are equivalent:
\begin{enumerate}
\item[(i) ]$\Sigma$ is an entire hyperplane.
\item[(ii) ]$\Sigma$ is a facet hyperplane of~$B$.
\item[(iii) ]There is no facet of~$\Sigma$ intersecting $J(\Sigma)$ in dimension $\dim(V)-2$.
\end{enumerate}
\end{lemma}

The following lemma is proved by a straightforward modification\footnote{Unfortunately, the statement of  \cite[Lemma~4.6]{sort_camb} is different enough that we must prove Lemma~\ref{half}, rather than simply quoting \cite[Lemma~4.6]{sort_camb}.}
 of the proof of  \cite[Lemma~4.6]{sort_camb}.

\begin{lemma}\label{half}
Let~$\B$ be the set of facet hyperplanes of the base region~$B$, and let $H\in\A$.
Then the following are equivalent:
\begin{enumerate}
\item[(i) ]$H$ contains exactly two shards.
\item[(ii) ]$H\not\in\B$ but there exists $\K\subseteq \B$ with $|\K|=2$ and $H\in\A_\K$.
\end{enumerate}
\end{lemma}
\begin{proof}
If (ii) holds then it is immediate from the definition that~$H$ is cut by the two elements of $\K$.
Both cuts remove the same subspace from~$H$, and Lemma~\ref{para cutting} implies that~$H$ is not cut by any other hyperplane.
Thus (i) holds.

Conversely, suppose (i) holds.
Then by Lemma~\ref{whole}, $H\not\in\B$.
Let $\K\subseteq\B$ be the set of hyperplanes in $\B$ which cut~$H$.
Suppose that $\K$ has exactly two hyperplanes $H'$ and $H''$.
Then by (i), these must cut~$H$ along the same codimension-2 subspace of~$V$, namely $H\cap H''$.
In this case,~$H$ is in $\A_{\set{H',H''}}$, and we have established (ii).
We complete the proof by showing that $|\K|=2$.

Each hyperplane $H'$ in $\K$ cuts~$H$ along some codimension-2 subspace~$U$ of~$V$ with $U\subseteq H'$.
By (i), this subspace~$U$ is the same for each $H'\in\K$.
In particular, the codimension-2 subspace~$U$ is contained in each $H'\in\K$.
Since~$B$ is a simplicial cone, the normal vectors to its facet-defining hyperplanes are linearly independent.
Thus $\bigcap\K$ has codimension~$|\K|$, and we conclude that $|\K|\le 2$.

Now suppose for the sake of contradiction that $|\K|<2$.
If $\K$ is a singleton, then let $H'$ be its unique element.
In this case $H'$ cuts~$H$.
Choose~$\Sigma$ to be the shard in~$H$ that is weakly below $H'$.
Thus the minimal element $J(\Sigma)$ of $U(\Sigma)$ is weakly below $H'$.
Let $H''\in(\B\setminus\K)$.
Since $H''$ does not cut~$H$, and since $H$ contains exactly two shards, there is some element of $U(\Sigma)$ that is weakly below $H''$.
Thus the minimal element $J(\Sigma)$ of $U(\Sigma)$ is weakly below $H''$.
But now $J(\Sigma)$ is a region weakly below every hyperplane in $\B$, so that $J(\Sigma)$ must be~$B$.
This is a contradiction, since~$B$ is not join-irreducible.

If $\K$ is empty, then let~$\Sigma$ be either of the two shards in~$H$.
Arguing as in the previous paragraph, we see that $J(\Sigma)$ is weakly below every hyperplane in $\B$ and reach the same contradiction.
\end{proof}

For any $H\in\A$, the \emph{depth} of~$H$ is the minimum cardinality of the separating set of a region separated from~$B$ by~$H$.
Suppose $J\in\R$ has $H\in S(J)$ and $|S(J)|=\depth(H)$.
If $J$ covers two or more other regions, at most one of those regions is separated from $J$ by~$H$, and thus it is possible to go down in the poset of regions while remaining separated from~$B$ by~$H$.
This contradition proves that any region $J$ with $H\in S(J)$ and $|S(J)|=\depth(H)$ must be join-irreducible in $(\R,\le)$.
Furthermore, $J$ is separated from $J_*$ by~$H$.
The following lemma makes possible an argument by induction on depth in the proof of Lemma~\ref{other direction}.

\begin{lemma}\label{depth}
If~$H$ is not a basic hyperplane of~$\A$ then there exists a rank-two full subarrangement $\A'$ containing~$H$ such that both basic hyperplanes of $\A'$ have depth strictly smaller than the depth of~$H$.
\end{lemma}
\begin{figure}
\scalebox{.65}{\includegraphics{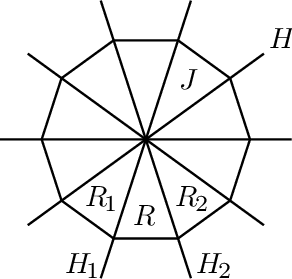}}
\caption{A figure for the proof of Lemma~\ref{depth}.}
\label{depthfig}
\end{figure}

\begin{proof}
Suppose~$H$ is not a basic hyperplane of~$\A$ and let $J$ be any region with $H\in S(J)$ and $|S(J)|=\depth(H)$.
Then Lemma~\ref{whole} says that there is a facet of $\Sigma(J)$ intersecting $J$ in dimension $\dim(V)-2$.
This intersection of $J$ with a facet of $\Sigma(J)$ is some codimension-2 face~$F$ of~$\F$.
The set of hyperplanes in~$\A$ containing~$F$ is a rank-two full subarrangement $\A'$.
Figure~\ref{depthfig} represents $\A'$ and the set of $\A$-regions containing~$F$.
Since the intersection $\bigcap\A'$ of the hyperplanes in $\A'$ contains a facet of $\Sigma(J)$, in particular~$H$ is not basic in $\A'$.
We claim that both basic hyperplanes~$H_1$ and~$H_2$ of $\A'$ have depth strictly less than the depth of~$H$.
Since the region $J$ contains~$F$, there is an $\A$-region~$R$ whose separating set (as an $\A$-region) is $S(J)\setminus\A'$.
The region~$R$ is covered by regions $R_1$ and $R_2$, also containing~$F$ and having respectively~$H_1\in S(R_1)$ and~$H_2\in S(R_2)$.
Since $J$ only covers one other region, and that cover is through~$H$ (not through~$H_1$ or~$H_2$), we have $J\not\in\set{R_1,R_2}$.
In particular, $|S(J)|>|S(R_1)|$.
But $\depth(H_1)\le|S(R_1)|$, so $\depth(H)=|S(J)|>\depth(H_1)$.
Similarly, $\depth(H)>\depth(H_2)$.
\end{proof}

The following two lemmas are the key technical ingredients in the proof of Proposition~\ref{recover cone}, which is crucial in the proofs in Section~\ref{int sec}.
They are roughly converse to each other.

\begin{lem}\label{other direction}
Let $\A'$ be a full subarrangement of~$\A$ and let $H\in(\A\setminus\A')$.
Let $\K$ be the set of basic hyperplanes of $\A'$ which are not cut by~$H$.
Let $H'\in(\A'\setminus(\A')_\K)$.
Then there exists a hyperplane $H''\in(\A\setminus\A')$ with $H''\cap(\bigcap\A')=H\cap(\bigcap\A')$ such that $H''$ cuts $H'$.
\end{lem}
\begin{proof}
We prove the lemma by reducing it to successively weaker statements.
First, we weaken the conclusion of the lemma by removing the requirement that $H''\cap(\bigcap\A')=H\cap(\bigcap\A')$.

\begin{weaker}\label{w1}
Let $\A'$ be a full subarrangement of~$\A$ and let $H\in(\A\setminus\A')$.
Let $\K$ be the set of basic hyperplanes of $\A'$ which are not cut by~$H$.
Let $H'\in(\A'\setminus(\A')_\K)$.
Then there exists a hyperplane $H''\in(\A\setminus\A')$ that cuts $H'$.
\end{weaker}

Given Weaker Assertion~\ref{w1}, the full lemma can be proved as follows:
Let $\A''$ be the smallest full subarrangement of~$\A$ containing $\A'$ and $H$.
This is the set of hyperplanes in~$\A$ which contain $H\cap(\bigcap\A)$.
By Weaker Assertion 1 (with $\A''$ replacing $\A$), there exists a hyperplane $H''\in(\A''\setminus\A')$ such that $H''$ cuts $H'$.
Then $H''\cap(\bigcap\A')=H\cap(\bigcap\A')=\bigcap\A''$. 

Next we strengthen the hypotheses by requiring that $H'$ is not contained in any proper standard subarrangement of $\A'$.
Assuming this additional hypothesis, the requirement that $H'\not\in(\A')_\K$ is equivalent to the requirement that some basic hyperplane of $\A'$ is cut by $H$.

\begin{weaker}\label{w2}
Let $\A'$ be a full subarrangement of~$\A$.
Let $H$ be a hyperplane in $(\A\setminus\A')$ which cuts some basic hyperplane of $\A'$.
Let $H'$ be a hyperplane in $\A'$ that is not contained in any proper standard subarrangement of $\A'$.
Then there exists a hyperplane $H''\in(\A\setminus\A')$ which cuts $H'$.
\end{weaker}

Given Weaker Assertion~\ref{w2}, we prove Weaker Assertion~\ref{w1} as follows.
Assuming the hypotheses of Weaker Assertion~\ref{w1}, let $(\A')_{\K'}$ be the smallest proper standard subarrangement of $\A'$ with $H'\in(\A')_{\K'}$.
Then $H'$ is not contained in any proper standard parabolic subgroup of $(\A')_{\K'}$.
The assumption that $H'\in(\A'\setminus(\A')_\K)$ implies that some basic hyperplane of $(\A')_{\K'}$ is cut by $H$.
Thus Weaker Assertion~\ref{w2} applies, with $\A'$ replaced by $(\A')_{\K'}$, and asserts that there exists a hyperplane $H''\in(\A\setminus(\A')_{\K'})$ such that $H''$ cuts $H'$.
Now Lemma~\ref{para cutting} says that $H''\in(\A\setminus\A')$.

Our final weakening of the statement specializes the hypotheses to a very special case: the case where the rank of~$\A$ is three and the rank of $\A'$ is two.
When $\A'$ has rank two, the hyperplane $H'$ is in a proper standard subarrangement of $\A'$ if and only if $H'$ is one of the two basic hyperplanes of $\A'$.

\begin{weaker}\label{w3}
Let~$\A$ be an arrangement of rank three and let $\A'$ be a full rank-two subarrangement of~$\A$.
Let $H$ be a hyperplane in $(\A\setminus\A')$ which cuts some basic hyperplane of $\A'$.
Let $H'$ be a non-basic hyperplane in $\A'$.
Then there exists a hyperplane $H''\in(\A\setminus\A')$ which cuts $H'$.
\end{weaker}

Given Weaker Assertion~\ref{w3}, we now prove Weaker Assertion~\ref{w2} by induction on the depth of $H'$ in $\A'$.
Assume the hypotheses of Weaker Assertion~\ref{w2} and let $d$ be the depth of $H'$ in $\A'$.
If $d=1$ then, since $H'$ is not contained in any proper standard parabolic subgroup of $\A'$, the rank of $\A'$ is one and $H'$ is the unique basic hyperplane of $\A'$.
Thus since $H$ cuts some basic hyperplane of $\A'$, Weaker Assertion~\ref{w2} holds with $H''=H$.

If $d>1$, then by Lemma~\ref{depth}, there is a full rank-two subarrangement $\widetilde{\A'}$ of $\A'$ containing $H'$, such that the basic hyperplanes~$H_1$ and~$H_2$ of $\widetilde{\A'}$ both have strictly smaller depth than $H'$.
Let $(\A')_{\K_1}$ be the smallest standard subarrangement containing~$H_1$ and let $(\A')_{\K_2}$ be the smallest standard subarrangement containing~$H_2$.
(Possibly $(\A')_{\K_1}=\A'$ or $(\A')_{\K_2}=\A'$ or both.)
The union $\K_1\cup\K_2$ must be $\B'$, the set of all basic hyperplanes of $\A'$;  otherwise,~$H_1$ and~$H_2$ are contained in the same proper standard subarrangement of $\A'$ and then $H'$ is also contained in the same proper standard subarrangement.
In particular, without loss of generality, $H$ cuts some basic hyperplane of $(\A')_{\K_1}$.
Also,~$H_1$ is not in any proper standard subarrangement of $(\A')_{\K_1}$, since $(\A')_{\K_1}$ is the smallest standard subarrangement of $\A'$ containing~$H_1$.

By induction on $d$, there exists a hyperplane $\widetilde{H}\in(\A\setminus(\A')_{\K_1})$ cutting~$H_1$.
By Lemma~\ref{para cutting}, $\widetilde{H}\in(\A\setminus\A')$, because no hyperplane in $\A'\setminus(\A')_{\K_1}$ cuts $H_1\in(\A')_{\K_1}$.
Consider the full subarrangement $\widetilde{\A}$ of~$\A$ consisting of hyperplanes containing $\widetilde{H}\cap(\bigcap\widetilde{\A'})$.
By Weaker Assertion~\ref{w3}, (with~$\A$, $\A'$ and~$H$ replaced by $\widetilde{\A}$, $\widetilde{\A'}$ and $\widetilde{H}$), there exists a hyperplane $H''\in(\widetilde{\A}\setminus\widetilde{\A'})$ which cuts $H'$.
If $H''\in\A'$, then the entire rank-three full subarrangement $\widetilde{\A}$ is contained in $\A'$, contradicting the fact that $\widetilde{H}\in\A\setminus\A'$.
Thus $H''\in\A\setminus\A'$.

We have shown that Weaker Assertion~\ref{w3} implies Weaker Assertion~\ref{w2}.
We complete the proof of the lemma by proving Weaker Assertion~\ref{w3}.
First, $\A'$ cannot be a standard subarrangement of $\A$, because if so, Lemma~\ref{para cutting} would imply that no hyperplane of $\A'$ is cut by $H$, contradicting the hypothesis that some basic hyperplane of $\A'$ is cut by $H$.
Now Lemma~\ref{half} implies that $H'$ is cut by some hyperplane $H''$ besides the basic hyperplanes of $\A'$.
Necessarily $H''\in(\A\setminus\A')$.
\end{proof}

\begin{lem}\label{one direction} 
Let $\A'$ be a full subarrangement of~$\A$ and let $H\in(\A\setminus\A')$.
Suppose $H$ cuts some hyperplane $H'\in\A'$.
Then there exists a hyperplane $H''\in(\A\setminus\A')$ with $H''\cap(\bigcap\A')=H\cap(\bigcap\A')$ such that $H''$ cuts some basic hyperplane of $\A'$.
\end{lem}
\begin{proof}
Exactly as in the proof of Lemma~\ref{other direction}, it is enough to prove the weaker assertion where the requirement $H''\cap(\bigcap\A')=H\cap(\bigcap\A')$ is removed from the conclusion.

Let $\B$ be the set of basic hyperplanes of~$\A$ and let $\B'$ be the set of basic hyperplanes of $\A'$.
We first claim that $\B'\not\subseteq\B$.
Supposing to the contrary that $\B'\subseteq\B$, the full subarrangement $\A'$ is the standard subarrangement $\A_{\B'}$.
But in this case, by Lemma~\ref{para cutting}, $H'$ is not cut by any hyperplane in $\A'\setminus(\A_{\B'})$.
This contradiction proves the claim.

The claim can be restated:  There exists a basic hyperplane $H'''$ of $\A'$ that is not basic in~$\A$.
By Lemma~\ref{whole}, $H'''$ is cut by some hyperplane $H''$ in~$\A$.
But since $H'''$ is basic in $\A'$, Lemma~\ref{whole} implies that $H''\in(\A\setminus\A')$.
\end{proof}

\section{Intersections of shards}\label{int sec}
In this section, we consider the set $\Psi(\A,B)$ of arbitrary intersections of shards and the lattice $(\Psi(\A,B),\supseteq)$, consisting of the elements of $\Psi(\A,B)$ partially ordered by reverse containment.
The space $V$ is in $\Psi(\A,B)$ by convention:  it is the intersection of the empty set of shards.
When it does not cause confusion, we write $\Psi$ instead of $\Psi(\A,B)$.
The main result of this section is that the set $\Psi$ is in bijection with the set~$\R$ of regions of~$\A$, so that the lattice $(\Psi,\supseteq)$ can be thought of as a partial order on~$\R$.
(Recall that, throughout, $\A$ is assumed to be simplicial.)

\begin{example}\label{i25Psi}
This is a continuation of Example~\ref{i25shards}.
When~$\A$ consists of five lines through the origin in $\reals^2$, the set $\Psi$ consists of ten cones, namely the origin, the eight shards shown in Figure~\ref{i25 fig}.c, and the whole space $\reals^2$.
\end{example}

\begin{figure}
\scalebox{.9}{\includegraphics{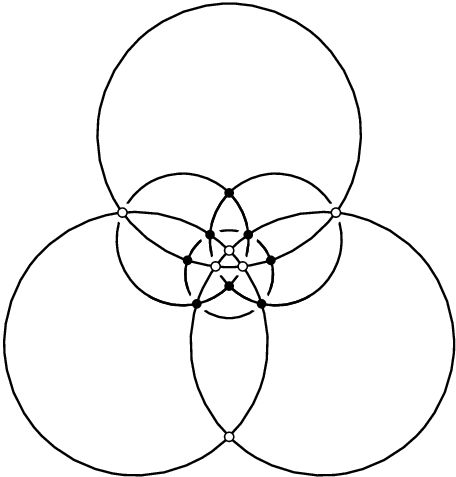}}
\caption{$\Psi$ in the example $\A=\A(S_4)$.}
\label{Psi A3}
\end{figure}

\begin{example}\label{Psi A3 ex}
This is a continuation of Examples~\ref{labeled S4} and~\ref{A3shards ex}.
When~$\A$ is the Coxeter arrangement $\A(S_4)$ in $\reals^3$, the elements of $\Psi$ are the origin, eleven one-dimensional cones (three of which are entire lines), the eleven shards shown in Figure~\ref{A3shards} and the whole space $\reals^3$.
Each cone intersects the unit sphere in one of six ways:  an empty intersection, a single point, a pair of antipodal points, an arc of a great circle, a great circle, or the entire sphere.
Figure~\ref{Psi A3} depicts these intersections in stereographic projection.
Thus the shards are shown as circles or circular arcs and the one-dimensional shards are pictured as points or pairs of points.
A white dot indicates a point which is paired with its antipodal point.
(To find antipodal points, note that any two of the circles shown intersect in a pair of antipodal points.)
\end{example}

Since each shard is a convex cone, the elements of $\Psi$ are all convex cones.
Each shard is a union of codimension-1 faces of the fan $\F=\F(\A)$.
Thus an intersection of shards is an intersection of unions of faces.
Distributing the intersection over the union, and keeping in mind that the intersection of shards is a convex cone, we have the following:
\begin{prop}\label{union of faces}
If $\Gamma\in\Psi$ has dimension $d$ then $\Gamma$ is a union of (closed) $d$-dimensional faces of the fan~$\F$.
\end{prop}

The key fact about shard intersections is the observation that a cone in $\Psi$ can be recovered from any of the full-dimensional faces it contains.
Given a face $F\in\F$, define a cone $\Gamma(F)\in\Psi$ as follows:
Let $[Q,R]$ be the interval in the poset of regions $(\R,\le)$ corresponding to~$F$.
Then $\Gamma(F)=\bigcap\set{\Sigma(P\covered R):Q\le P\covered R}$.  

\begin{prop}\label{recover cone}
If $F\in\F$ and $\Gamma\in\Psi$ have $\dim(F)=\dim(\Gamma)$ and $F\subseteq \Gamma$ then $\Gamma=\Gamma(F)$.
Furthermore $\Gamma$ is the intersection of all shards containing~$F$.
\end{prop}

\begin{proof}
By Proposition~\ref{facial}, each shard in $\set{\Sigma(P\covered R):Q\le P\covered R}$ contains a different facet of~$R$, and the number of covers in $\set{(P\covered R):Q\le P\covered R}$ is the codimension of~$F$.
Thus $\Gamma(F)$ contains~$F$ and has the same dimension as~$F$.
Furthermore, the subspace $U=\bigcap\set{H(P\covered R):Q\le P\covered R}$ is the smallest subspace containing~$F$.
Since~$F$ is a full dimensional subset of $\Gamma$,~$U$ is also the smallest subspace containing $\Gamma$.

Proposition~\ref{cov ref basic} states that the hyperplanes $\set{H(P\covered R): P\covered R}$ are the basic hyperplanes of a full subarrangement.
The set $\set{H(P\covered R):Q\le P\covered R}$ is weakly smaller than $\set{H(P\covered R): P\covered R}$, so $\set{H(P\covered R):Q\le P\covered R}$ is the set of basic hyperplanes of a weakly smaller full subarrangement $\A'$ consisting of all hyperplanes containing~$U$.
Let $\B'$ be this set of basic hyperplanes of $\A'$.

Since $\Gamma\in\Psi$, we can write $\Gamma=\bigcap\set{\Sigma_1,\ldots,\Sigma_k}$ for some shards $\Sigma_i$.
Alternately, the cone $\Gamma$ is obtained as follows:
We cut~$U$ along every hyperplane not containing~$U$ that cuts any hyperplane $H(\Sigma_i)$ for the defining shards $\Sigma_i$.
Each of the resulting pieces is a union of faces of~$\F$.
The piece containing~$F$ is $\Gamma$.
Similarly, $\Gamma(F)=\bigcap\set{\Sigma(P\covered R):Q\le P\covered R}$ is obtained from~$U$ by cutting~$U$ along every hyperplane not containing~$U$ that cuts any of the hyperplanes in $\B'$.
Again, the piece containing~$F$ is $\Gamma(F)$.
To prove that $\Gamma(F)=\Gamma$, we show that both of these cutting schemes cut~$U$ in exactly the same way.

On the one hand, suppose there exists a hyperplane $H\in\A\setminus\A'$ cutting some $H(\Sigma_i)$.
By Lemma~\ref{one direction}, there exists a hyperplane $H'\in\A\setminus \A'$ which cuts a hyperplane in $\B'$, with $H\cap U=H'\cap U$.
On the other hand, suppose a hyperplane $H\in\A\setminus\A'$ cuts some hyperplane in $\B'$.
Then the intersection of the basic hyperplanes of $\A'$ which are not cut by~$H$ is strictly larger than~$U$.
In particular, since $U=\bigcap\set{H(\Sigma_1),\ldots,H(\Sigma_k)}$, there is some hyperplane $H(\Sigma_i)$ which does not contain the intersection of the basic hyperplanes of $\A'$ which are not cut by~$H$.
By Lemma~\ref{other direction}, with $H'=H(\Sigma_i)$, there exists a hyperplane $H''$ cutting $H(\Sigma_i)$ such that $H''\cap U=H\cap U$.

We have proved the first assertion of the proposition.
Now, let $\Gamma'$ be the intersection of all shards containing~$F$.
Then $F\subseteq \Gamma'\subseteq \Gamma$, so~$F$ is a full-dimensional face contained in $\Gamma'$, and by the first statement of the proposition, $\Gamma'=\Gamma(F)=\Gamma$.
\end{proof}

The number of covers in $\set{(P\covered R):Q\le P\covered R}$ is the codimension of~$F$.
Thus if $\Gamma$ has codimension~$d$ then Proposition~\ref{recover cone} expresses $\Gamma$ as the intersection of $d$ distinct shards.
Call these shards the \emph{canonical shards} containing $\Gamma$.
Note that the choice of canonical shards containing $\Gamma$ is well-defined:
Any choice of~$F$ in Proposition~\ref{recover cone} yields a set of $\codim(\Gamma)$-many shards contained in the basic hyperplanes of full subarrangement $\A'=\set{H\in\A:\Gamma\subseteq H}$.
There is a unique such set of shards whose intersection is $\Gamma$.

\begin{prop}\label{shard face}
If $\Gamma\in\Psi$ then any face of $\Gamma$ is in $\Psi$.
\end{prop}
\begin{proof}
We begin by showing that any facet of a shard~$\Sigma$ is in $\Psi$.
The facet~$C$ is the intersection of~$\Sigma$ with some hyperplane~$H_1$ that cuts the hyperplane $H'=H(\Sigma)$.
Then~$H_1$ is a basic hyperplane of the rank-two full subarrangement $\A'$ containing $H'$ and~$H_1$.
Let~$H_2$ be the other basic hyperplane of $\A'$.
Since~$H_1$ cuts $H'$, we know that~$H_2\neq H'$.

Let~$F$ be a face of~$\F$ with $F\subseteq C$ and $\dim(F)=\dim(C)$.
The intersection $\Gamma(F)$ of all shards containing~$F$ is contained in~$\Sigma$, because~$\Sigma$ contains~$F$.
Since $\Gamma(F)$ is a convex polyhedral cone contained in~$\Sigma$ and the linear span of $\Gamma(F)$ equals the linear span of the face~$C$ of~$\Sigma$, we conclude that $\Gamma(F)\subseteq C$.

Proposition~\ref{recover cone} implies that the shard intersections contained in~$C$ are the pieces obtained by cutting the subspace $\bigcap\A'$ along all hyperplanes that cut either $H_1$ or $H_2$.
Thus if $\Gamma(F)$ is properly contained in~$C$ then there exists a hyperplane $H\in(\A\setminus\A')$ which cuts either $H_1$ or $H_2$ and intersects the relative interior of~$C$.
Then Lemma~\ref{other direction} states that there exists a hyperplane $H''\in(\A\setminus\A')$ with $H''\cap(\bigcap\A')=H\cap(\bigcap\A')$ such that $H''$ cuts $H'$.
But then $H''$ intersects the relative interior of~$C$ as well and, since $H'=H(\Sigma)$, $H''$ intersects the relative interior of~$\Sigma$.
This contradicts that fact that~$\Sigma$ is a single shard,
By this contradiction, we conclude that the containment $\Gamma(F)\subseteq C$ is in fact equality.
In particular $C\in\Psi$.

Next we observe that any facet of a cone $\Gamma\in\Psi$ is in $\Psi$.
Write $\Gamma=\bigcap\set{\Sigma_1,\ldots,\Sigma_k}$ where $\Sigma_1,\ldots,\Sigma_k$ are the canonical shards containing $\Gamma$.
Then $\Gamma$ is the subset of the subspace $\bigcap\set{H_{\Sigma_1},\ldots,H_{\Sigma_k}}$ defined by all of the facet-defining inequalities of all the shards $\Sigma_i$.
In particular, a facet~$C$ of $\Gamma$ is defined by some facet-defining inequality for a facet of some $\Sigma_i$.
Thus, by the special case already proved,~$C$ is the intersection of $\Gamma$ with some shard intersection, so that $F\in\Psi$.

Finally, if~$F$ is a lower-dimensional face of $\Gamma$, it is the intersection of a set of facets of $\Gamma$, and thus $F\in\Psi$ as well.
\end{proof}

\begin{prop}\label{standard face}
Let $\Gamma\in\Psi$.
Then there is a unique minimal subset $\K$ of $\B$ such that $\Gamma$ is an intersection of shards contained in hyperplanes in $\A_\K$.
The minimal face of $\Gamma$ is $\bigcap\K$.
\end{prop}
\begin{proof}
If $\Gamma$ can be expressed as an intersection of shards contained in hyperplanes in $\A_\K$, then Lemma~\ref{para cutting} says that $\bigcap\K\subseteq\Gamma$.
Since $\Gamma$ is a convex cone, there is a unique minimal $\K\subseteq\B$ such that $\bigcap\K\subseteq\Gamma$.
Then every shard containing $\Gamma$ also contains $\bigcap\K$, so $H(\Sigma)$ contains $\bigcap\K$, or in other words, $H(\Sigma)\in\A_\K$.
Thus the minimal $\K$ with $\bigcap\K\subseteq\Gamma$ is the desired subset.

To show that $\bigcap\K$ is the minimal face of $\Gamma$, it is enough to consider the special case where~$\Sigma$ is a shard.
Indeed, given the special case, an intersection $\Gamma=\Sigma_1\cap\cdots\cap\Sigma_k$ of shards has as its minimal face some intersection of subspaces $(\bigcap\A_{\K_1})\cap\cdots(\bigcap\A_{\K_k})$ where $\A_{\K_i}$ is the minimal subset of $\B$ such~$\Sigma$ is contained in a hyperplane in $\A_{\K_i}$.
Then $(\bigcap\A_{\K_1})\cap\cdots(\bigcap\A_{\K_k})=\bigcap\A_{\K}$, where $\K=\K_1\cup\cdots\cup\K_k$ is the minimal subset of $\B$ such $\Gamma$ is an intersection of shards contained in hyperplanes in $\A_\K$.

We now reduce the special case to an even more special case.
Let~$\Sigma$ and $\K$ be as in the special case.
Lemma~\ref{para cutting} says that $\Sigma\supseteq(\bigcap\K)$.
Lemma~\ref{para cutting} says, furthermore, that~$\Sigma$ is obtained by cutting a hyperplane containing $\bigcap\K$ along its intersections with other hyperplanes containing $\bigcap\K$.
By Lemma~\ref{full cutting}, we can ignore the rest of $\A$, or in other words, reduce to the case where $\A_\K=\A$, or equivalently $\K=\B$.
Thus we need only show that, when $H(\Sigma)$ is not containined in any proper standard subarrangement, $\bigcap\A$ is a face of~$\Sigma$.

Suppose for the sake of contradiction that $\bigcap\A$ is not a face of~$\Sigma$.
The minimal face of any polyhedral cone is a subspace.
Let~$U$ be the minimal face of~$\Sigma$ and let $\A'$ be the full subarrangement consisting of hyperplanes in~$\A$ containing~$U$.
Since $U\supsetneq(\bigcap\A)$, we have $\A'\subsetneq\A$.
Since~$U$ is the minimal face, both $H(\Sigma)$ and every hyperplane cutting $H(\Sigma)$ is in $\A'$.
Because $H(\Sigma)$ is not contained in any proper standard subarrangement, at least one basic hyperplane of $\A'$ is not basic in~$\A$.
Thus by Lemma~\ref{whole}, there is a hyperplane $H$ in $\A\setminus\A'$ that cuts some basic hyperplane of $\A'$.
Now Lemma~\ref{other direction} (or, more conveniently, Weaker Assertion~\ref{w2} in the proof of Lemma~\ref{other direction}) implies that some hyperplane in $\A\setminus\A'$ cuts $H(\Sigma)$.
In particular, $U\not\subseteq\Sigma$, and this contradiction shows that $\bigcap\A$ is a face of~$\Sigma$.
\end{proof}

The \emph{lattice of shard intersections} is the set $\Psi$ partially ordered by reverse containment.
The unique minimal element of $(\Psi,\supseteq)$ is the empty intersection, interpreted as the entire space~$V$.
The unique maximal element is the intersection of the set of all shards.
This maximal element coincides with $\bigcap\A$.
The poset $(\Psi,\supseteq)$ is a lattice.
The join operation is intersection, and since the poset also has a unique minimal element, meets can be defined in the usual way in terms of joins:
$\Gamma_1\meet\Gamma_2$ is the intersection of all shards~$\Sigma$ such that $\Gamma_1\subseteq\Sigma$ and $\Gamma_2\subseteq\Sigma$.

We now show that $(\Psi,\supseteq)$ is in fact a partial order on the set~$\R$ of regions of~$\A$, by giving an explicit bijection between~$\R$ and $\Psi$.
Define a map $\psi:\R\to\Psi$ sending~$R$ to the intersection of the lower shards of~$R$.
In light of Theorem~\ref{R can join rep}, $\psi(R)$ is the intersection of all shards $\Sigma(J)$ such that $J$ is a canonical joinand of~$R$ in $(\R,\le)$.
Define a map $\rho:\Psi\to\R$ by setting $\rho(\Gamma)=\Join_{\Sigma\supseteq \Gamma}J(\Sigma)$, with the join taken in the poset of regions $(\R,\le)$.

\begin{prop}\label{bijection}
\noindent
\begin{enumerate}
\item[(i) ]$\psi$ is a bijection from~$\R$ to $\Psi$ with inverse map $\rho$.
\item[(ii) ]$\rho$ is an order-preserving map from $(\Psi,\supseteq)$ to the poset of regions $(\R,\le)$.
\item[(iii) ]The number of lower hyperplanes of $R\in\R$ equals the codimension of $\psi(R)$.
\end{enumerate}
\end{prop}

\begin{proof}
We first show that $\psi:\R\to\Psi$ is surjective.
Let $\Gamma$ be a $d$-dimensional cone in $\Psi$, and let~$F$ be a $d$-dimensional face of~$\F$ contained in $\Gamma$.
Then by Proposition~\ref{recover cone}, $\Gamma=\bigcap\set{\Sigma(P\covered R):Q\le P\covered R}$ for some $Q$ and~$R$ in~$\R$.
The canonical join-representation of~$R$ is $\set{J(P\covered R):P\covered R}=\set{J(\Sigma(P\covered R):P\covered R}$.
Thus the smaller set $\set{J(\Sigma(P\covered R)):Q\le P\covered R}$ is the canonical join-representation of some element $R'$.
We have $\Gamma=\psi(R')$.

We next show that $\rho\circ\psi$ is the identity map on~$\R$.
Given $R\in\R$, we have $R=\Join_{J\in\Can(R)}J$ in $(\R,\le)$.
Thus since $\psi(R)\subseteq\Sigma(J)$ for each $J\in\Can(R)$, it is enough to show that for any shard~$\Sigma$ containing $\psi(R)$, we have $J(\Sigma)\le R$.
But for each $J\in\Can(R)$, the shard $\Sigma(J)$ contains a facet of the region~$R$.
Thus $\psi(R)$ contains the face~$F$ of~$R$ obtained by intersecting the facets of~$R$ separating~$R$ from regions $Q$ having $Q\covered R$.
By Proposition~\ref{facial}, the set $\set{P\in\R:P\supseteq F}$ is an interval $I$ in $(\R,\le)$, and~$R$ is the maximal element of $I$.
Now, given any shard~$\Sigma$ containing $\psi(R)$, in particular~$\Sigma$ contains~$F$, so there is a region $P\in U(\Sigma)$ such that $P$ contains~$F$.
Thus $P\in I$ so that in particular $P\le R$.
But since $J(\Sigma)$ is the unique minimal element of $U(\Sigma)$, we conclude that $J(\Sigma)\le R$.

Since $\psi$ is surjective and $\rho\circ\psi$ is the identity map on~$\R$, the map $\psi$ is a bijection from~$\R$ to $\Psi$ with inverse map $\rho$.
This is (i).

Next we establish (ii).
If $\Gamma_1\supseteq \Gamma_2$ then in particular, the set of shards~$\Sigma$ such that $\Sigma\supseteq \Gamma_1$ is contained in the set of shards~$\Sigma$ such that $\Sigma\supseteq \Gamma_2$ and therefore $\rho(\Gamma_1)\le\rho(\Gamma_2)$.

The codimension of $\psi(R)$ is the size of the canonical join representation of~$R$, which equals the number of lower hyperplanes of~$R$ by Theorem~\ref{R can join rep}.
This is (iii).
\end{proof}

\begin{example}\label{Psi A3 bij ex}
This is a continuation of Examples~\ref{labeled S4}, \ref{A3shards ex}, and~\ref{Psi A3 ex}.
In Figure~\ref{Psi A3}, each region~$R$ appears as a curvilinear triangle.
The cone $\psi(R)$ is the intersection of the shards containing the concave edges of the triangle.
\end{example}

The lattice induced on~$\R$, via the bijection of Proposition~\ref{bijection}, from $(\Psi(\A,B),\supseteq)$ is denoted by $(\R,\preceq)$ or $(\R(\A),\preceq_B)$.
We call $(\R,\preceq)$ the \emph{shard intersection order} on~$\R$.
The unique minimal element of $(\R,\preceq)$ is~$B$, and the unique maximal element is $-B$.
We emphasize that $\preceq$ corresponds to $\supseteq$, not to $\subseteq$.
We also emphasize that the partial orders $(\R,\le)$ and $(\R,\preceq)$ are distinct.
By Proposition~\ref{bijection}(ii), the shard intersection order $(\R,\preceq)$ is a weaker order than the poset of regions $(\R,\le)$.

\begin{example}
This is a continuation of Examples~\ref{i25shards} and~\ref{i25Psi}.
(See Figure~\ref{i25 fig}.c.)
When~$\A$ consists of five lines in $\reals^2$, the poset $(\Psi,\supseteq)$ has $\reals^2$ as its unique minimal element and the origin as its unique maximal element.
The $8$ ($1$-dimensional) shards are pairwise incomparable under containment, and occur at rank 1 (i.e.\ codimension~1).
Thus the poset $(\R(\A),\preceq)$ has~$B$ as its unique minimal element and $-B$ as its unique maximal element.
The other $8$ regions~$\R$ are pairwise incomparable and occur at at rank 1.
\end{example}

\begin{figure}
\scalebox{.8}{\includegraphics{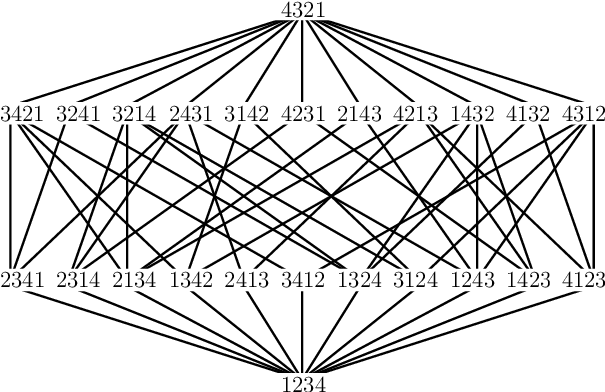}}
\caption{$(S_4,\preceq)$.}
\label{S4 preceq}
\end{figure}

\begin{example}
Continuing Examples~\ref{labeled S4}, \ref{A3shards ex}, \ref{Psi A3 ex} and~\ref{Psi A3 bij ex}, the lattice $(\R,\preceq)$ is shown in Figure~\ref{S4 preceq} for the Coxeter arrangement $\A(S_4)$.
The elements of $(\R,\preceq)$ are labeled by the corresponding permutations, as shown in Figure~\ref{A3labels}.
\end{example}

\section{Properties of the shard intersection order}\label{prop sec}
In this section, we establish some basic properties of the shard intersection order for a simplicial hyperplane arrangement.
It will be convenient to pass freely between $(\R,\preceq)$ and $(\Psi,\supseteq)$.

\begin{prop}\label{graded}
The lattice $(\R,\preceq)$ is graded, with the rank of $R\in\R$ equal to the number of lower hyperplanes of~$R$. 
Alternately, the rank of $\Gamma\in\Psi$ is the codimension of $\Gamma$.
\end{prop}
\begin{proof}
The minimal element~$V$ of $(\Psi,\supseteq)$ has codimension zero.
Suppose $\Gamma\supseteq \Gamma'$ in $(\Psi,\supseteq)$.
Let $F'$ be some full-dimensional face in $\Gamma'$.
By Lemma~\ref{tech}, there is a face~$F$ of~$\F$ which is full-dimensional in $\Gamma$, such that $F'$ is a face of~$F$.
If $\dim(\Gamma)=\dim(\Gamma')$, then $F=F'$, so $\Gamma=\Gamma(F')=\Gamma'$ by Proposition~\ref{recover cone}.
If $\dim(\Gamma)>\dim(\Gamma')+1$ then $\dim(F)>\dim(F')+1$, so there is a face $G$ of~$F$ with $F'\subsetneq G\subsetneq F$, and $\Gamma(G)$ is an element of $\Psi$ with $\Gamma'\subsetneq \Gamma(G)\subsetneq \Gamma$.
This proves the first assertion, and the second assertion follows by Proposition~\ref{bijection}(iii).
\end{proof}

\begin{prop}\label{atomic coatomic}
The lattice $(\R,\preceq)$ is atomic and coatomic.
\end{prop}
\begin{proof}
It is immediate that an element of $(\Psi,\supseteq)$ is join-irreducible if and only if it is a shard.
Since all of the shards are atoms, $(\Psi,\supseteq)\cong(\R,\preceq)$ is an atomic lattice.

Let $k=\dim(\bigcap\A)$.
We will show that every element $\Gamma$ of $(\Psi,\supseteq)$ of dimension at least $k+2$ contains at least two distinct elements of $\Psi$ whose dimension is $\dim(\Gamma)-1$.
By Proposition~\ref{graded}, this implies that the only meet-irreducible elements of $(\Psi,\supseteq)$ are the coatoms.

Let $\Gamma\in\Psi$ have dimension at least $k+2$.
If $\Gamma$ has two or more facets, then we are done by Propositions~\ref{shard face} and~\ref{graded}.
If $\Gamma$ has no facets, then $\Gamma$ has no proper faces, so by Proposition~\ref{standard face}, $\Gamma$ is the subspace $\bigcap\K$ for some $\K\subseteq\B$.
Furthermore $|\K|\le|\B|-2$, because $\Gamma\in\Psi$ has dimension at least $k+2$.
Let $H_1$ and $H_2$ be distinct hyperplanes in $\B\setminus\K$.
Both $H_1$ and $H_2$ are shards by Lemma~\ref{whole}.
Thus $\Gamma\cap H_1$ and $\Gamma\cap H_2$ are distinct shard intersections of dimension $\dim(\Gamma)-1$.
If $\Gamma$ has one facet, then $\Gamma$ has exactly one proper face~$F$.
By Proposition~\ref{standard face},~$F$ is the subspace $\bigcap\K$ for some $\K\subseteq\B$ with $|\K|\le|\B|-1$.
Let $H$ be a hyperplane in $\B\setminus\K$.
Then~$F$ and $\Gamma\cap H$ are distinct shard intersections of dimension $\dim(\Gamma)-1$.
\end{proof}

We omit the easy proof of the following proposition.
\begin{prop}\label{reducible}
If~$\A$ is the direct sum of $\A_1$ and $\A_2$ then $(\R(\A),\preceq_B)$ is isomorphic to $(\R(\A_1),\preceq_{B_1})\times(\R(\A_2),\preceq_{B_2})$, where $B_i$ is the $\A_i$-region containing~$B$, for $i=1,2$. 
\end{prop}
As a special case of Proposition~\ref{reducible}, if~$W$ is a reducible finite Coxeter group with $W\cong W_1\times W_2$, then the shard intersection order on~$W$ is isomorphic to the product of the shard intersection orders on $W_1$ and $W_2$.

For each $R\in\R$, define $L(R)=\Meet\set{P:P\covered R}$.
This is the maximal element in $(\R,\le)$ which is below~$R$ but which does not contain any lower hyperplane of~$R$ in its separating set.

\begin{prop}\label{auto}
The map $\Gamma\mapsto -\Gamma$ is an automorphism of $(\Psi,\supseteq)$.
\end{prop}
\begin{proof}
The operation of cutting a hyperplane into shards has antipodal symmetry, so that $\Sigma\mapsto-\Sigma$ is an involution on the set of shards.
Thus $\Gamma\mapsto -\Gamma$ is an involution on $\Psi$.
The map is also containment-preserving, so it is an automorphism.
\end{proof}

We now explain the relationship between the shard intersection poset $(\Psi,\supseteq)$ and two other geometrically-defined lattices associated to $\A$: the intersection lattice $(\Int(\A),\supseteq)$ of~$\A$ and the face lattice $(\F(\A),\subseteq)$ of the fan $\F(\A)$ associated to~$\A$.
For convenience, we think of the whole space~$V$ as a face of the fan~$\F$.
(Alternately, we may work with the zonotope that is dual to~$\F$ and take the usual convention that the empty set is a face of any polytope.)

Any order preserving map $\eta:P\rightarrow Q$ defines a relation on the set $\bar{P}=\set{\eta^{-1}(q):q\in Q}$ of fibers of $\eta$ as follows:
Set $F_1 \le_{\bar{P}}F_2$ if there exist $a \in F_1$ and $b \in F_2$ such that $a \le_P b$.
If $\le_{\bar{P}}$ is a partial order on $\bar{P}$, then $(\bar{P},\le_{\bar{P}})$ is called the {\em fiber poset}\footnote{Often this poset is called the \emph{quotient} of $P$ with respect to $\eta$.
However, because lattice-theoretic quotients play a prominent role in what follows, we prefer the term fiber poset.} of the map $\eta:P\rightarrow Q$.

Given a face $F\in\F$, recall that $\Gamma(F)$ is the intersection of all shards containing~$F$.
Given a set $X\subseteq V$, let $U(X)$ be the subspace obtained as the intersection of all hyperplanes in~$\A$ containing $X$.

\begin{prop}\label{psi fib}
The intersection lattice $(\Int(\A),\supseteq)$ is isomorphic to the fiber poset of $U:\Psi\to \Int(\A)$.
The lattice $(\Psi,\supseteq)$ is anti-isomorphic to the fiber poset of $\Gamma:\F\to\Psi$.
\end{prop}

The proof of Proposition~\ref{psi fib} is simple.
We borrow some terminology from \cite{recursion}.
An \emph{order projection} is an order-preserving map $\eta:P\rightarrow Q$, with the following property: 
For all $x \le y$ in $Q$, there exist $a \le b \in P$ with $\eta(a)=x$ and $\eta(b)=y$.
In particular, an order-projection is surjective.

When an order preserving map $\eta:P\rightarrow Q$ has the property that $\le_{\bar{P}}$ is a partial order, there is a surjective order-preserving map $\nu:P \rightarrow \bar{P}$ given by $\nu:a \mapsto \eta^{-1}(\eta(a))$, and an injective order-preserving map $\bar{\eta}:\bar{P} \rightarrow Q$ such that $\eta = \bar{\eta}\circ \nu$.
The following easy fact about order projections appears as~\cite[Proposition~1.1]{recursion}.
\begin{prop}
\label{op properties}
Let $\eta:P\rightarrow Q$ be an order projection.
Then 
\begin{enumerate}
\item[(i) ] the relation $\le_{\bar{P}}$ is a partial order, and 
\item[(ii) ] $\bar{\eta}$ is an isomorphism of posets.
\end{enumerate}
\end{prop}

\begin{proof}[Proof of Proposition~\ref{psi fib}]
The map~$U$ is order-preserving and we show that it is an order projection.
Given $U_1\supseteq U_2$ in $\Int(\A)$, choose a face $G$ of~$\F$ contained in $U_2$ having full dimension in $U_2$.
By Lemma~\ref{tech}, there is a face~$F$ with $G\subseteq F\subseteq U_1$ such that~$F$ is full-dimensional in $U_1$.
Then $\Gamma(F)\supseteq \Gamma(G)$, $U(\Gamma(F))=U_1$ and $U(\Gamma(G))=U_2$.
Thus~$U$ is an order projection, so we apply Proposition~\ref{op properties} to obtain the first assertion.

The proof of the second assertion is very similar, except that the map $\Gamma$ is order-reversing.
Let $\Gamma_1\supseteq \Gamma_2$ in $\Psi$, choose a face $G$ of~$\F$ contained in~$\Gamma_2$ having full dimension in~$\Gamma_2$.
Since $\Gamma_1\supseteq \Gamma_2$, Lemma~\ref{tech} says that there is a face~$F$ with $G\subseteq F\subseteq \Gamma_1$ such that~$F$ is full-dimensional in~$\Gamma_1$.
Then by Proposition~\ref{recover cone}, $\Gamma(F)=\Gamma_1$ and $\Gamma(G)=\Gamma_2$.
Thus $\Gamma$ is an order projection from~$\F$ to the dual of $(\Psi,\supseteq)$, so we apply Proposition~\ref{op properties} to see that the fiber poset of $\Gamma:\F\to\Psi$ is isomorphic to $(\Psi,\subseteq)$, the dual of $(\Psi,\supseteq)$.
\end{proof}

For any $Q,R\in\R$ with $Q\le R$, let $I(Q,R)=\set{J(P\covered P'):Q\le P\covered P'\le R}$.

\begin{prop}\label{preceq combin}
Let $Q,R\in\R$.  Then $Q\preceq R$ if and only if $I(L(Q),Q)\subseteq I(L(R),R)$.
\end{prop}
\begin{proof}
For any $R\in\R$, the interval $[L(R),R]$ coincides with $\set{P\in\R:P\supseteq F}$, where~$F$ is the intersection of~$R$ with all regions $R'$ such that $R'\covered R$.
The cone $\psi(R)$ contains~$F$ and $\dim(F)=\dim(\psi(R))$.
Thus by Proposition~\ref{recover cone}, a shard~$\Sigma$ contains $\psi(R)$ if and only if it contains~$F$.
But~$\Sigma$ contains~$F$ if and only if it separates two adjacent cones $P$ and $P'$ with $P,P'\in [L(R),R]$.
Thus the sets $I(L(R),R)$ and $\set{J(\Sigma):\psi(R)\subseteq\Sigma}$ coincide.

By definition, $Q\preceq R$ if and only if $\psi(Q)\supseteq\psi(R)$.
This holds if and only if the set of shards $\set{\Sigma:\psi(Q)\subseteq\Sigma}$ is contained in the set of shards $\set{\Sigma:\psi(R)\subseteq\Sigma}$.
The latter occurs if and only if $I(L(Q),Q)\subseteq I(L(R),R)$.
\end{proof}

\begin{prop}\label{lower int}
If $R\in\R$, then the lower interval $[1,R]$ in $(\R(\A),\preceq_B)$ is isomorphic to $(\R(\A'),\preceq_{B'})$, where $\A'$ is the full subarrangement of~$\A$ consisting of hyperplanes containing $\bigcap\Lower(R)$, and $B'$ is the $\A'$-region containing~$B$.
\end{prop}

\begin{proof}
The lower interval $[1,R]$ corresponds to the lower interval $[V,\psi(R)]$ in $(\Psi,\supseteq)$.
The face $G=R\cap\bigcap\Lower(R)$ of~$\F$ is contained in $\psi(R)$ and is full-dimensional in $\psi(R)$.
By Proposition~\ref{recover cone}, $\psi(R)$ is the intersection of all shards containing $G$, so a cone $\Gamma\in\Psi$ is in $[V,\psi(R)]$ if any only if it contains $G$.
In light of Lemma~\ref{full cutting}, taking $x$ to be a point in the relative interior of $G$, the shards of $(\A',B')$ coincide, in a small neighborhood of $x$, with the shards of $(\A,B)$ containing $G$.
This coincidence defines a bijection $\Gamma\mapsto \Gamma'$ from the interval $[V,\psi(w)]$ in $(\Psi(\A,B),\supseteq)$ to the poset $(\Psi(\A',B'),\supseteq)$.
We show that this bijection is an isomorphism of posets.
Let $\Gamma_1,\Gamma_2\in\Psi(\A,B)$ have $\Gamma_1\supseteq G$ and $\Gamma_2\supseteq G$.

Suppose $\Gamma_1\supseteq \Gamma_2$.
We use Proposition~\ref{union of faces} and Lemma~\ref{tech} to find faces $F_1,F_2\in\F$ with $F_1\supseteq F_2\supseteq G$ such that $F_1$ is contained full-dimensionally in~$\Gamma_1$ and $F_2$ is contained full-dimensionally in~$\Gamma_2$.
Furthermore,~$\Gamma_1$ is the intersection of all shards containing $F_1$ and~$\Gamma_2$ is the intersection of all shards containing $F_2$.
Let $F_1'$ be the face of $\F(\A')$ such that $F_1\subseteq F_1'$ and $\dim(F_1')=\dim(F_1)$.
Define $F_2'$ similarly.
By Proposition~\ref{recover cone}, $\Gamma_1'$ is the intersection of all shards of $\A'$ containing $F_1'$ and $\Gamma_2'$ is the intersection of all shards of $\A'$ containing $F_2'$.
Since $F_1'\supseteq F_2'$, we have $\Gamma_1'\supseteq \Gamma_2'$.

Conversely, suppose $\Gamma_1'\supseteq \Gamma_2'$.
We find faces $F_1',F_2'\in\F'$ with $F_1'\supseteq F_2'$ such that~$F_1'$ is contained full-dimensionally in $\Gamma_1'$ and $F_2'$ is contained full-dimensionally in~$\Gamma_2'$.
By Proposition~\ref{recover cone}, $\Gamma_1'$ is the intersection of all shards of $(\A',B')$ containing~$F_1'$ and $\Gamma_2'$ is the intersection of all shards of $(\A',B')$ containing $F_2'$.
Let $F_1$ be the face of $\F(\A)$ such that $G\subseteq F_1\subseteq F_1'$ and $\dim(F_1)=\dim(F_1')$.
Define $F_2$ similarly.
By Proposition~\ref{recover cone},~$\Gamma_1$ is the intersection of all shards of~$\A$ containing $F_1$ and~$\Gamma_2$ is the intersection of all shards of~$\A$ containing $F_2$.
All the shards containing~$F_1$ or $F_2$ also contain $G$.
Since $F_1\supseteq F_2$, we have $\Gamma_1\supseteq \Gamma_2$.
\end{proof}

Proposition~\ref{W lower int} is a more detailed version of Proposition~\ref{lower int} in the Coxeter case.
\begin{proof}[Proof of Proposition~\ref{W lower int}]
For any $w\in W$, let $W'$ be the parabolic subgroup generated by the set $S'=\cov(w)$.
Proposition~\ref{lower int} says that $[1,w]_\preceq$ is isomorphic to the lattice $(W',\preceq)$ defined with respect to the Coxeter system $(W',S')$.
The Coxeter system $(W',S')$ is isomorphic to $(W_J,J)$ and Proposition~\ref{W lower int} follows.
\end{proof}

Recall from Section~\ref{arr sec} that the parabolic subset $\R_\K$ is the set of regions whose separating sets are contained in $\A_K$.

\begin{theorem}\label{mobius}
The M\"{o}bius number of $(\R,\preceq)$ is $(-1)^{\rank(\A)}$ times the number of regions in~$\R$ that are not contained in any proper parabolic subset of~$\R$.
Equivalently, by inclusion-exclusion,
\[\mu_\preceq(B,-B)=\sum_{\K\subseteq \B}(-1)^{|\K|}\left|\R_\K\right|.\]
\end{theorem}

\begin{proof}
Let $\nu(B,-B)$ be the proposed M\"{o}bius number.
By Lemma~\ref{RK}, we rewrite
\[\nu(B,-B)=\sum_{\K\subseteq \B}(-1)^{|\K|}\,\bigl|\bigl\lbrace Q\in\R:Q\supseteq \bigl((-B)\cap\bigcap\K\bigr)\bigr\rbrace\bigr|.\]

For each $R\in\R$, let $G(R)$ be the intersection of all facets of~$R$ separating~$R$ from regions covered by~$R$.
Let $\A'(R)$ be the full subarrangement of~$\A$ consisting of hyperplanes $H\in\A$ with $H\supseteq G(R)$.
By Proposition~\ref{cov ref basic}, $\A'(R)$ has basic hyperplanes $\Lower(R)$.
By Proposition~\ref{facial}, the $\A'(R)$-regions are in one-to-one correspondence with $\A$-regions containing $G(R)$.
Thus by Proposition~\ref{lower int}, we have
\[\nu(B,R)=\sum_{\K\subseteq \Lower(R)}(-1)^{|\K|}\,\bigl|\bigl\lbrace Q\in\R:Q\supseteq \bigl(R\cap\bigcap\K\bigr)\bigr\rbrace\bigr|.\]
The map $\K\mapsto\bigl(R\cap\bigcap\K\bigr)$ is a bijection between subsets of $\Lower(R)$ and faces of~$R$ containing $G(R)$.
Thus $\nu(B,R)$ is
\[\sum_{\substack{F\in\F\\R\supseteq F\supseteq G(R)}}(-1)^{\codim(F)}\,\bigl|\bigl\lbrace Q\in\R:Q\supseteq F\bigr\rbrace\bigr|=\sum_{\substack{F\in\F\\R\supseteq F\supseteq G(R)}}(-1)^{\codim(F)}\sum_{\substack{Q\in\R\\Q\supseteq F}}1.\]

By Theorem~\ref{shelling}, any linear extension of the dual $(\R,\ge)$ of the poset of regions is a shelling order on the simplicial fan~$\F$.
A standard result uses this shelling order to obtain a partition of the faces of~$\F$ into intervals in the face poset of~$\F$.
The intervals obtained are exactly the intervals $[G(R),R]$, so for any face~$F$ of~$\F$, there exists a unique region~$R$ such that $R\supseteq F\supseteq G(R)$.
Thus, to show that $\sum_{R\in\R}\nu(B,R)$ vanishes, we reverse the order of summation and obtain
\[\sum_{Q\in\R}\,\sum_{\substack{F\in\F\\F\subseteq Q}}(-1)^{\codim(F)}\,\sum_{\substack{R\in\R\\R\supseteq F\supseteq G(R)}}1\,\,=\,\,
\sum_{Q\in\R}\,\sum_{\substack{F\in\F\\F\subseteq Q}}(-1)^{\codim(F)}\]
The inner sum is zero because the face lattice of the cone $Q$ is Eulerian, with rank function given by dimension.
\end{proof}

Theorem~\ref{W mobius} is a special case of Theorem~\ref{mobius}.
However, there is an alternate proof of Theorem~\ref{W mobius} which is interesting in comparison to the proof of Theorem~\ref{mobius}:
While the proof of Theorem~\ref{mobius} relies on the fact that face lattices of polytopes are Eulerian posets, the proof below rests on inclusion/exclusion.
In other words, it rests on the fact that the Boolean lattice is Eulerian.
\begin{proof}[Alternate proof of Theorem~\ref{W mobius}]
In light of Proposition~\ref{W lower int}, we must verify that the following sum vanishes:
\[\sum_{w\in W}\sum_{J\subseteq \Des(w)}(-1)^{|J|}\left|W_J\right|=\sum_{J\subseteq S}(-1)^{|J|}\left|W_J\right|\sum_{\substack{w\in W\\J\subseteq \Des(w)}}1.\]
The inner sum is the number of maximal-length representatives of cosets of $W_J$ in~$W$. 
This number is $|W|/\left|W_J\right|$, so the sum reduces to zero.
\end{proof}

Proposition~\ref{W lower int} can also be used to give a computationally effective recursive formula for counting maximal chains in $(W,\preceq)$, namely Proposition~\ref{num max}.

\begin{proof}[Proof of Proposition~\ref{num max}]
The number of maximal chains in $(W,\preceq)$ is the sum over all coatoms $w$ of $(W,\preceq)$ of the number of maximal chains in $[1,w]$.
Every coatom $w$ is a maximal-length coset representative of the subgroup $W_\br{s}$ for some unique $s\in S$.
On the other hand, for each $s\in S$, every coset of $W_\br{s}$ has a unique maximal-length coset representative.
This representative $w$ has $\rank(W)-1$ descents and thus is a coatom of $(W,\preceq)$, except if $w$ is $w_0$, which has $\rank(W)$ descents.
For each $s\in S$, there are exactly $\left|W\right|/\left|W_\br{s}\right|$ cosets of $W_\br{s}$, and exactly one of these cosets has $w_0$ as its maximal length representative.
The proposition now follows by Proposition~\ref{W lower int}.
\end{proof}

The fan $\F(\A)$ defined by a central hyperplane arrangement~$\A$ is dual to a zonotope.
In fact, there is an uncountable family of zonotopes dual to $\F(\A)$.  
We fix some particular dual zonotope to call $Z(\A)$.
We now describe the connection between the order complex of $(\R,\preceq)$ and a certain pulling triangulation of $Z(\A)$.
This connection is inspired by~\cite{Loday}, in a way that is easier to explain in Section~\ref{nc sec}.

Given a total order of the vertices on some polytope, the \emph{pulling triangulation} of the polytope is defined recursively as follows:
Let $v_0$ be the first vertex in the total order.
Recursively triangulate each face~$F$ not containing the vertex $v_0$, by pulling, using the restriction of the total order to vertices in~$F$.
This triangulated polyhedral complex is extended to a triangulation of the entire polytope by coning at the vertex $v_0$.
In the pulling triangulation, every face of the polytope, even a face containing~$F$, is triangulated by pulling, using the restriction of the total order to vertices in~$F$.
(See~\cite{Lee} for more details.)

More generally, it is enough to partially order the vertices of the polytope, as long as each face of the polytope has a unique minimal vertex in the partial order.
We define the triangulation $\Delta(\A)$ of $Z(\A)$ to be the pulling triangulation with respect to the dual poset of regions $(\R,\ge)$.

We define a map $\delta$ from chains in $(\R,\preceq)$ to simplices in $\Delta(\A)$.
Let $\chi$ be a chain in $(\R,\preceq)$.
Place the vertex $v_0$ in $\delta(\chi)$ if and only if $-B\in \chi$, and let $\chi'=\chi\setminus\set{-B}$.
If $\chi'=\emptyset$ then $\delta(\chi)$ is either $\emptyset$ or $\set{v_0}$.
Otherwise, let~$R$ be the maximal element of $\chi'$ and let~$F$ be the unique maximal face of $Z(\A)$ such that~$R$ is dual to the top vertex of~$F$ in $(\R,\le)$.
By~Proposition~\ref{lower int}, there is a bijection between elements of $[B,R]_\preceq$ and vertices of~$F$.
We inductively map the chain $\chi'$ to a simplex in the triangulation of~$F$.
This process is easily reversible. Thus:
\begin{theorem}\label{zon chain tri}
The map $\delta$ is a dimension-preserving bijection between the order complex of $(\R,\preceq)$ and the pulling triangulation $\Delta(\A)$.
\end{theorem}
In particular, the two complexes have the same $f$-vector.

\begin{figure}
\includegraphics{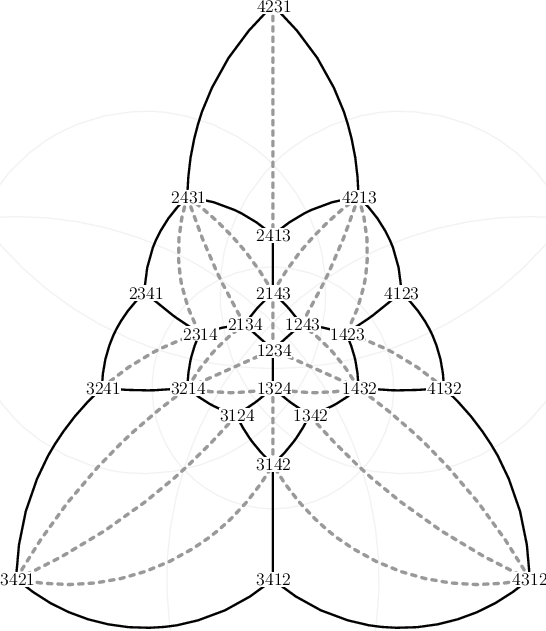}
\caption{A stage in the construction of $\Delta(\A(S_4))$.}
\label{A3 pull fig}
\end{figure}

\begin{example}\label{perm tri example}
Figure~\ref{A3 pull fig} depicts a stage in the construction of $\Delta(\A)$ in the case $\A=\A(S_4)$, in the same stereographic projection employed in previous figures.
The labels and solid black lines show the $1$-skeleton of the polyhedral subcomplex obtained by removing all faces containing the vertex $v_0$ (labeled by $4321$) from $Z(\A)$.
The dotted gray lines illustrate the recursive triangulation of the subcomplex.
Combinatorially, $\Delta(\A)$ is the cone over the two-dimensional simplicial complex shown.
Theorem~\ref{zon chain tri} says that the 34 maximal chains of $(\Psi,\preceq)$ correspond to (cones over the) 34 triangles shown in the figure.
\end{example}

\section{Shards and lattice congruences}\label{cong sec}
In this section, we provide background information about congruences of a finite lattice.
In particular, we describe the connection between shards and lattice congruences of the poset of regions $(\R,\le)$.
(As mentioned in Section~\ref{arr sec}, the poset of regions is a lattice because of our assumption that $\A$ is simplicial.)
We emphasize that we are studying lattice congruences on the poset of regions $(\R,\le)$, not on the shard intersection order $(\R,\preceq)$.
Indeed, it appears that in general the shard intersection order has no interesting lattice congruences.

A congruence on a lattice~$L$ is an equivalence relation on~$L$ such that $x_1\equiv x_2$ and $y_1\equiv y_2$ implies that $(x_1\meet y_1)\equiv(x_2\meet y_2)$ and $(x_1\join y_1)\equiv(x_2\join y_2)$.
We use the symbol~$\Theta$ to represent a typical lattice congruence, and write $[x]_\Theta$ for the $\Theta$-class of $x\in L$.
The poset $\Con(L)$ of all congruence relations on a finite lattice~$L$, partially ordered by refinement, is known (see e.g. \cite[Theorem~II.3.11]{Gratzer}) to be a distributive lattice.

When~$L$ is finite, congruences on~$L$ have an order-theoretic characterization which is easily verified, or which follows from more general results in~\cite{Cha-Sn}.

\begin{prop}\label{cong char}
An equivalence relation~$\Theta$ on a finite lattice~$L$ is a congruence if and only if the following three conditions hold:
\begin{enumerate}
\item[(i) ]Each $\Theta$-class is an interval in~$L$.
\item[(ii) ]The downward projection may $\pidown^\Theta$, sending $x\in L$ to the bottom element of $[x]_\Theta$, is order-preserving.
\item[(iii) ]The upward projection may $\piup_\Theta$, sending $x\in L$ to the top element of $[x]_\Theta$, is order-preserving.
\end{enumerate}
\end{prop}

The lattice quotient $L/\Theta$ is the lattice whose elements are the $\Theta$-classes, with $[x]_\Theta\join[y]_\Theta=[x\join y]_\Theta$ and $[x]_\Theta\meet[y]_\Theta=[x\meet y]_\Theta$.
Equivalently, $L/\Theta$ is the partially ordered set whose elements are the $\Theta$-classes, with $[x]_\Theta\le[y]_\Theta$ if and only if there exists $x'\in[x]_\Theta$ and $y'\in[y]_\Theta$ with $x'\le y'$.

Three additional facts will be useful.
The first and second are known and easily verified, and the third is \cite[Prop 2.2]{con_app}.
\begin{prop}\label{bottoms}
For any congruence~$\Theta$ on a finite lattice~$L$, the lattice quotient $L/\Theta$ is isomorphic, as a partially ordered set, to the induced subposet $\pidown^\Theta(L)$ of~$L$.
\end{prop}

\begin{prop}\label{bottom char}
Let~$L$ be a finite lattice and let $x\in L$ have a canonical join representation.
Let~$\Theta$ be a lattice congruence on~$L$, with associated downward projection $\pidown^\Theta$.
Then $x$ is in $\pidown^\Theta(L)$ if and only if no canonical joinand of $x$ is contracted.
\end{prop}


\begin{prop}\label{cover lemma}
Let~$L$ be a finite lattice, let~$\Theta$ be a congruence on~$L$, and let $x\in L$.
Then the map $y\mapsto[y]_\Theta$ restricts to a one-to-one correspondence between elements of~$L$ covered by $\pidown^\Theta(x)$ and elements of $L/\Theta$ covered by $[x]_\Theta$.
\end{prop}

When an edge $x\covered y$ in the Hasse diagram of $L$ has $x\equiv y$, we say that the edge is \emph{contracted}.
Proposition~\ref{cong char} says, in particular, that when $L$ is finite, congruence classes are intervals.
Thus we can completely decribe a congruence~$\Theta$ on a finite lattice~$L$ by listing the edges that are contracted by~$\Theta$.

\begin{figure}
\scalebox{.75}{\includegraphics{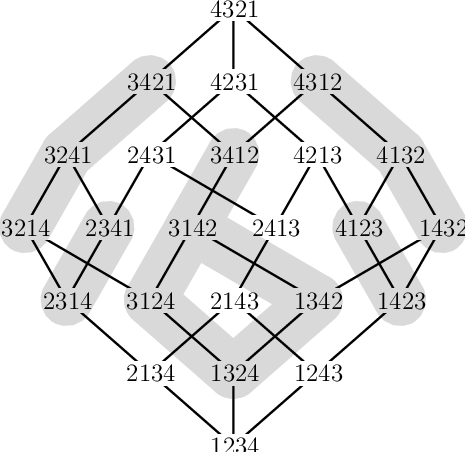}}
\caption{A Cambrian congruence on the weak order on $S_4$.}
\label{A3camb_cong}
\end{figure}

\begin{example}\label{cong ex}
Figure~\ref{A3camb_cong} shows a congruence on the weak order on $S_4$.
The congruence relation is depicted by shading all contracted edges.
Thus the unshaded vertices are singleton congruence classes, and the shading groups the remaining vertices into congruence classes.
This is an example of a Cambrian congruence, as we explain in Example~\ref{S4 camb}.
\end{example}

Edges cannot be contracted independently. 
Rather, contracting one edge may require other edges to be contracted in order to obtain a congruence.
Recall that Proposition~\ref{shard ji} describes a bijection between shards in~$\A$ and join-irreducible regions in $(\R,\le)$.
For each shard~$\Sigma$, the corresponding join-irreducible region $J(\Sigma)$ is the unique minimal region (in the sense of the poset of regions) among upper regions of~$\Sigma$.
For each join-irreducible region $J$, the corresponding shard $\Sigma(J)$ is the shard containing the intersection of $J$ and $J_*$, where $J_*$ is the unique region covered by $J$ in $(\R,\le)$.
When a congruence~$\Theta$ contracts $J_*\covered J$, we say that~$\Theta$ contracts $J$.

\begin{prop}\label{shard contract}
Let $Q\covered R$ in $(\R,\le)$.
Then a lattice congruence~$\Theta$ on $(\R,\le)$ contracts $Q\covered R$ if and only if it contracts $J(\Sigma(Q\covered R))$.
Thus if $\Sigma(Q'\covered R')=\Sigma(Q\covered R)$ then~$\Theta$ contracts $Q\covered R$ if and only if~$\Theta$ contracts $Q'\covered R'$.
\end{prop}
\begin{proof}
The key points are the following:
\begin{enumerate}
\item $H(Q\covered R)$ is in $S(J)$ and in $S(R)$;
\item $S(J_*)=S(J)\setminus\set{H(Q\covered R)}$;
\item $S(Q)=S(R)\setminus\set{H(Q\covered R)}$;
\item $J\le R$; and 
\item $J_*\le Q$.
\end{enumerate}
The first three points are immediate from the hypotheses, the fourth point follows by Proposition~\ref{shard ji} and the last point follows from the first four.

These five points enable the following simple argument for the first assertion of the proposition.
If $Q\equiv R$ modulo~$\Theta$ then since~$\Theta$ is a lattice congruence, we have $J_*=Q\meet J\equiv R\meet J=J$.
Conversely, if $J\equiv J_*$ then $R=J\join Q\equiv J_*\join Q=Q$.
The second assertion is immediate from the first.
\end{proof}

Proposition~\ref{shard contract} says that, given a shard~$\Sigma$ and a congruence~$\Theta$ on $(\R,\le)$ either all of the edges associated to~$\Sigma$ are contracted or none of the edges associated to~$\Sigma$ are contracted.
When the edges associated to~$\Sigma$ are contracted, we say that~$\Sigma$ is \emph{removed} by~$\Theta$.
(The reason we speak of ``removing'' shards rather than ``contracting'' shards will become clear below, when we make the connection between lattice congruences and fans.)

Proposition~\ref{shard contract} also makes it clear that a lattice congruence~$\Theta$ on $(\R,\le)$ is completely determined by the set of shards~$\Theta$ removes.
Equivalently,~$\Theta$ is completely determined by the set of join-irreducible elements it contracts.
The latter fact is easily proved for general finite lattices.
(See for example \cite[Theorem~2.30]{FreeLattices}.)
The key, then, to characterizing lattice congruences is to determine which sets of join-irreducible elements can be the set of join-irreducible elements contracted by a congruence.

\begin{example}\label{i2 cong}
Consider a rank-two hyperplane arrangement~$\A$ with $k$ hyperplanes and choose a base region~$B$.
The case $k=5$ was considered in Examples~\ref{labeled i25}, \ref{i25shards}, and~\ref{i25Psi}.
Label the regions~$B$, $-B$, $Q_1$, $Q_2$, \ldots, $Q_{k-1}$, $R_1$, $R_2$, \ldots, and $R_{k-1}$, as in Figure~\ref{i25 fig}.
Each of these regions is join-irreducible in $(\R,\le)$ except $\pm B$.

One can verify that, for any $i$ from 2 to $k-1$, there is a congruence contracting $Q_i$ and no other join-irreducible element.
In fact, the congruence contracts no other edge.
The same is true for $R_i$ with $i$ from 2 to $k-1$.
On the other hand, suppose some congruence~$\Theta$ contracts $Q_1$.
(That is, it sets $Q_1\equiv (Q_1)_*=B$.)
Then $-B=(Q_1\join R_1)\equiv (B\join R_1)=R_1$, and we conclude that~$\Theta$ contracts all $R_i$ with $i$ from 2 to $k-1$.
Furthermore, $Q_{k-1}=(-B\meet Q_{k-1})\equiv(R_1\meet Q_{k-1})=B$, and we conclude that~$\Theta$ contracts all $Q_i$ with $i$ from 2 to $k-1$.
Thus contracting $Q_1$ forces all $Q_i$ and $R_i$ with $i$ from 2 to $k-1$ to be contracted, and one can verify that $R_1$ is not forced to be contracted.
Similarly, contracting $R_1$ forces the same set of additional join-irreducibles to be contracted.

We verify Proposition~\ref{shard contract} in this example.
One of the two nontrivial assertions of the proposition in this example is that a congruence sets $B\equiv Q_1$ if and only if it sets $R_{k-1}\equiv -B$. 
The ``only if'' direction was already verified above, and the ``if'' direction is verified by a dual argument.
The other nontrivial assertion, that $B\equiv R_1$ if and only if it sets $Q_{k-1}\equiv -B$, is proved similarly.
\end{example}

For congruences on $(\R,\le)$, at least some of the forcing relations among shards have a nice geometric description.
Given two shards~$\Sigma$ and $\Sigma'$, say $\Sigma\to\Sigma'$ if $H(\Sigma)$ cuts $H(\Sigma')$ and $\Sigma\cap\Sigma'$ has codimension~2.
The digraph thus defined on shards is called the \emph{shard digraph}.
The following theorem is a restatement of part of \cite[Theorem~25]{hyperplane}.
We remind the reader that throughout the paper, the arrangement~$\A$ is assumed to be simplicial.

\begin{theorem}\label{acyclic shard cong}
Suppose the shard digraph defined by $(\A,B)$ is acyclic and let $X$ be a set of shards.
Then there exists a lattice congruence of $(\R,\le)$ contracting the shards in $X$ and no other shards if and only if $X$ is an order ideal in the transitive closure of the shard digraph.
\end{theorem}

\begin{figure}
\psfrag{2134}[cc][cc]{\large$\Sigma(2134)$}
\psfrag{1324}[cc][cc]{\large$\Sigma(1324)$}
\psfrag{1243}[cc][cc]{\large$\Sigma(1243)$}
\psfrag{2314}[cc][cc]{\large$\Sigma(2314)$}
\psfrag{3124}[cc][cc]{\large$\Sigma(3124)$}
\psfrag{1342}[cc][cc]{\large$\Sigma(1342)$}
\psfrag{1423}[cc][cc]{\large$\Sigma(1423)$}
\psfrag{2341}[cc][cc]{\large$\Sigma(2341)$}
\psfrag{2413}[cc][cc]{\large$\Sigma(2413)$}
\psfrag{3412}[cc][cc]{\large$\Sigma(3412)$}
\psfrag{4123}[cc][cc]{\large$\Sigma(4123)$}
\scalebox{.8}{\includegraphics{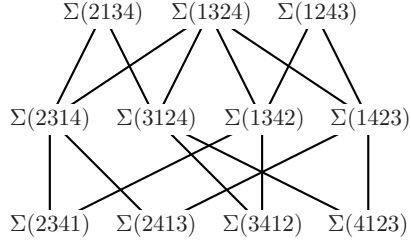}}
\caption{The shard digraph for $\A(S_4)$.}
\label{A3digraph fig}
\end{figure}

\begin{example}\label{A3sharddigraph}
Continuing Examples~\ref{A3shards ex} and~\ref{A3 shard ji}, the transitive closure of the shard digraph of the Coxeter arrangement $\A(S_4)$ is shown in Figure~\ref{A3digraph fig}.
(Cf. \cite[Figure~4]{congruence}.)
Each cover relation shown in Figure~\ref{A3digraph fig} is, of course, an arrow in the shard digraph (with arrows pointing down).
In addition, the shard digraph has an arrow from $\Sigma(2134)$ to each of the four shards at the bottom rank of the poset, and an arrow from $\Sigma(1243)$ to each of the same four shards.
\end{example}

The shard digraph of $(\A,B)$ can have directed cycles (see \cite[Figure~5]{hyperplane}), so Theorem~\ref{acyclic shard cong} does not apply in general.
In the motivating case of Coxeter arrangements, the shard digraph is always acyclic \cite[Proposition~28]{hyperplane}, so Theorem~\ref{acyclic shard cong} applies.
When the shard digraph has cycles, the usual construction produces a poset on a set of equivalence classes of vertices in the digraph.
One naturally wonders whether, in general, congruences of $(\R,\le)$ correspond to order ideal in the latter poset.
At present, we do not have a complete answer to this question.
However, the following fact suffices for the purposes of this paper.

\begin{prop}\label{good enough}
Let~$\Theta$ be a congruence on $(\R,\le)$ and let $X$ be the set of shards removed by~$\Theta$.
If~$\Sigma$ and $\Sigma'$ are shards with $\Sigma\in X$ and $\Sigma\to\Sigma'$, then $\Sigma'\in X$.
\end{prop}

\begin{proof}
Suppose $\Sigma\in X$ and $\Sigma\to\Sigma'$, and let~$F$ be a codimension-2 face of~$\F$ contained in $\Sigma\cap\Sigma'$.
Then Proposition~\ref{facial} says that the set of all regions in~$\R$ containing~$F$ is some interval $[P,P']$ in $(\R,\le)$.
By definition of the shard digraph, $H(\Sigma)$ is a lower hyperplane of~$P'$ and $H(\Sigma')$ is not.
The interval $[P,P']$ is isomorphic to a poset of regions of rank two.
Then $P$ and $P'$ correspond to~$B$ and $-B$ in Example~\ref{i2 cong}.
Let the additional regions in $[P,P']$ be labeled $Q_i$ and $R_i$ as in Example~\ref{i2 cong}.
Then one of the regions $Q_i$ or $R_i$ with $i$ between 2 and $k-1$ is an upper region of $\Sigma'$, separated by $\Sigma'$ from $Q_{i-1}$ or $R_{i-1}$.
Since $\Sigma\in X$, the restriction of~$\Theta$ to $[P,P']$ contracts (without loss of generality) $P'\covered Q_1$.
As explained in Example~\ref{i2 cong},~$\Theta$ must contract all edges of the form $Q_{i-1}\covered Q_{i}$ and $R_{i-1}\covered R_{i}$ with $i$ between 2 and $k-1$.
Thus by Proposition~\ref{shard contract},~$\Theta$ removes $\Sigma'$.
\end{proof}

For each $H\in\B$, let $R(H)$ stand for the region whose separating set is $\set{H}$.
The regions $R(H)$ are the atoms of the poset of regions $(\R,\le)$.

\begin{theorem}
\label{parabolic quotient}
Let~$\B$ be the set of basic hyperplanes of~$\A$.
Let $\K\subseteq\B$ and let $B'$ be the $\A_\K$-region containing~$B$.
Then the map $R\mapsto R_\K$ is a lattice homomorphism from $(\R(A),\le_B)$ to $(\R(\A_\K),\le_{B'})$.
The fibers of this homomorphism constitute the finest lattice congruence of $(\R(\A),\le_B)$ with $B\equiv R(H)$ for every $H\in(\B-\K)$.
This congruence contracts exactly those shards which lie in hyperplanes not in $\A_\K$.
\end{theorem}

The first assertion of Theorem~\ref{parabolic quotient} is \cite[Proposition~6.3]{congruence}. 
(Cf. \cite[Lemmas~4.3,~4.5]{Jed}.)
The second assertion is \cite[Theorem~6.9]{congruence}.
The third assertion follows immediately from the definition of $R_\K$.

Each lattice congruence~$\Theta$ on $(\R,\le)$ defines a complete fan $\F/\Theta$, refined by~$\F$.
The definition of $\F/\Theta$, and all of the properties listed, are quoted from~\cite[Sections~4--5]{con_app}.
By definition, two maximal cones of~$\F$ are in the same maximal cone of $\F/\Theta$ if and only if they are congruent mod~$\Theta$.  
In particular, $(\R,\le)/\Theta$ is a partial order on the maximal cones of $\F/\Theta$.
For any face $F\in(\F/\Theta)$, the set of maximal faces of $\F/\Theta$ containing~$F$ is an interval $[X,Y]$ in $(\R,\le)/\Theta$.
Any linear extension of $(\R,\le)/\Theta$ (or of its dual) is a shelling order on the maximal cones of $\F/\Theta$.
There exists a regular CW-sphere $\S(\F/\Theta)$ which is dual to $\F/\Theta$, in the sense that the two face posets are anti-isomorphic.
The Hasse diagram of the lattice quotient $(\R,\le)/\Theta$ is an orientation of the $1$-skeleton of $\S(\F/\Theta)$.

Let $R\in\R$ and let~$C$ be the maximal cone of $\F/\Theta$ containing~$R$.
Then the facet-defining hyperplanes of~$C$ consist of the lower hyperplanes of $\pidown^\Theta(R)$, which separate~$C$ from maximal cones covered by~$C$ in $(\R,\le)/\Theta$, and the upper hyperplanes of $\piup_\Theta(R)$, which separate~$C$ from maximal cones covering~$C$ in $(\R,\le)/\Theta$.

Shards are not mentioned in~\cite{con_app}, but the definition of $\F/\Theta$ has a straightforward rephrasing in terms of shards:
The maximal cones of $\F/\Theta$ are the closures of the connected components of $V\setminus\bigcap\set{\Sigma:\Sigma\mbox{ is not removed by $\Theta$}}$.
\section{Lattice congruences and shard intersections}\label{cong psi sec}
In this section, we define a shard intersection poset $(\pidown^\Theta(\R),\preceq)$ for each lattice congruence~$\Theta$ on $(\R,\le)$.
Furthermore, we show that the properties of $(\R,\preceq)$ are inherited by $(\pidown^\Theta(\R),\preceq)$.
The motivating example, where~$\A$ is a Coxeter arrangement,~$\Theta$ is a Cambrian congruence, and $(\pidown^\Theta(\R),\preceq)$ is the noncrossing partition lattice, is discussed in Section~\ref{nc sec}.
The poset $(\pidown^\Theta(\R),\preceq)$ is the restriction of $(\R,\preceq)$ to $\pidown^\Theta(\R)$, the set of bottom elements of $\Theta$-classes.
We approach $(\pidown^\Theta(\R),\preceq)$ via a geometrically defined partial order, as we approached $(\R,\preceq)$.
Define $\Psi/\Theta$ to be the set of intersections of shards \textbf{not} removed by~$\Theta$.

\begin{figure}
\scalebox{.9}{\includegraphics{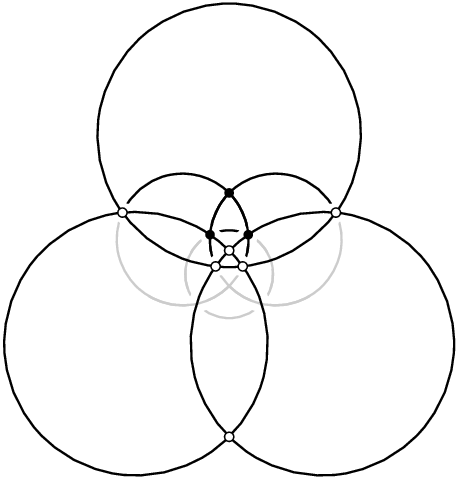}}
\caption{$\Psi/\Theta$ for the Cambrian congruence~$\Theta$ of Figure~\ref{A3camb_cong}.}
\label{Psi S4 bip}
\end{figure}

\begin{example}\label{Psi A3 bip example}
Figure~\ref{Psi S4 bip} depicts $\Psi/\Theta$, where~$\A$ is the Coxeter arrangement $\A(S_4)$ and~$\Theta$ is the congruence of the weak order on $S_4$ pictured in Figure~\ref{A3camb_cong}.
The removed shards are faded but not completely gone, to allow easy comparison with Figure~\ref{Psi A3}. 
\end{example}

We now show that the key properties of shard intersections (discussed in Section~\ref{int sec}), carry over to intersections of unremoved shards.
We also show that $(\pidown^\Theta(\R),\preceq)\cong(\Psi/\Theta,\supseteq)$.
Results about $(\pidown^\Theta(\R),\preceq)$ then follow by simple modifications of the proofs in Section~\ref{prop sec}.

\begin{prop}\label{Theta union of faces}
If $\Gamma\in\Psi/\Theta$ has dimension $d$ then $\Gamma$ is a union of (closed) $d$-dimensional faces of the fan $\F/\Theta$.
\end{prop}
\begin{proof}
It is enough to show that each unremoved shard~$\Sigma$ is a union of closed codimension~1 faces of $\F/\Theta$.
The general statement then follows just as as argued for Proposition~\ref{union of faces}.

Suppose for the sake of contradiction that some shard~$\Sigma$ and some \mbox{codimension-1} face~$F$ of $\F/\Theta$ have an intersection that is full-dimensional in~$\Sigma$ and in~$F$, but $F\not\subseteq\Sigma$.
Then some facet~$C$ of~$\Sigma$ intersects the relative interior of~$F$, and thus there is a hyperplane $H'$ cutting $H(\Sigma)$ such that $H'$ intersects the relative interior of~$F$.
Now~$F$ is a union of faces of~$\F$, since $\F$ refines $\F/\Theta$, and~$C$ is a union of faces of~$\F$ by Proposition~\ref{shard face} and Proposition~\ref{union of faces}. 
Thus their intersection is a union of faces of~$\F$, so we can choose a face $G$ of~$\F$ (not of $\F/\Theta$) of codimension~2 contained in~$C$ and~$F$.
This face $G$ is also contained in $H'$.
Because $H'$ is not cut along its intersection with $H(\Sigma)$ and since each shard in $H'$ is a union of faces of~$\F$, there is a unique shard $\Sigma'$ in $H'$ containing $G$.
We have $\Sigma'\to\Sigma$ in the shard digraph, so by Proposition~\ref{good enough}, $\Sigma'$ is not removed by~$\Theta$.
There is some codimension-1 face of~$\F$ containing $G$ and contained in $H'$, and since~$\F$ is refined by $\F/\Theta$, there is some codimension-1 face of $\F/\Theta$ containing $G$ and contained in $H'$.
But this face intersects the relative interior of~$F$, contradicting the fact that $\F/\Theta$ is a fan.
\end{proof}

Consider a cover relation $X\covered Y$ in $(\R,\le)/\Theta$.
Since the Hasse diagram of $(\R,\le)/\Theta$ is an orientation of the dual sphere to $\F/\Theta$, $X$ and $Y$ are adjacent maximal cones of $\F/\Theta$.
Their intersection is a codimension-1 face of $\F/\Theta$, so Proposition~\ref{Theta union of faces} implies that $X\cap Y$ is contained in some shard, which we represent by the symbol $\Sigma_\Theta(X\covered Y)$.

\begin{lemma}\label{canon canon}
Let $\Gamma\in\Psi/\Theta$ and let~$F$ be a face in $\F/\Theta$ with $F\subseteq\Gamma$ and $\dim(F)=\dim(\Gamma)$.
Let $[X,Z]$ be the interval in $(\R,\le)/\Theta$ corresponding to~$F$.
Then the set of canonical shards of $\Gamma$ is $\set{\Sigma_\Theta(Y\covered Z):X\le Y\covered Z}$.
\end{lemma}
\begin{proof}
Let~$R$ be the minimal region of the $\Theta$-class represented by $Z$.
By Proposition~\ref{cover lemma}, the map $\eta$ taking a region $Q$ covered by~$R$ to the maximal cone of $\F/\Theta$ containing $Q$ is a bijection between elements covered by~$R$ in $(\R,\le)$ and maximal cones of $\F/\Theta$ covered by $Z$ in $(\R,\le)/\Theta$.
The intersection of~$R$ and $Q$ is contained in the intersection of $Z$ and $\eta(Q)$.
Thus $\set{\Sigma_\Theta(Y\covered Z):X\le Y\covered Z}=\set{\Sigma(P\covered R):Q\le P\covered R}$.
\end{proof}

For $F\in\F/\Theta$, define $\Gamma_\Theta(F)=\bigcap\set{\Sigma_\Theta(Y\covered Z):X\le Y\covered Z}$, where $[X,Z]$ is the interval in $(\R,\le)/\Theta$ corresponding to $F$.
\begin{prop}\label{Theta recover cone}
If $F\in\F/\Theta$ and $\Gamma\in\Psi/\Theta$ have $\dim(F)=\dim(\Gamma)$ and \mbox{$F\subseteq \Gamma$} then $\Gamma=\Gamma_\Theta(F)$.
Furthermore $\Gamma$ is the intersection of all unremoved shards containing~$F$.
\end{prop}
\begin{proof}
Since~$\F$ refines $\F/\Theta$, there is a face $G$ of~$\F$ with $G\subseteq F$ and $\dim(G)=\dim(F)$.
Then by Proposition~\ref{recover cone}, $\Gamma=\Gamma(G)$, the intersection of the canonical shards of $\Gamma$.
By Lemma~\ref{canon canon}, $\Gamma=\Gamma_\Theta(F)$.
The second assertion follows exactly as in the proof of Proposition~\ref{recover cone}.
\end{proof}

\begin{prop}\label{canon unremoved}
Let $\Gamma\in\Psi$ and let~$\Theta$ be a lattice congruence on $(\R,\le)$.
Then~$\Gamma$ is in $\Psi/\Theta$ if and only if none of its canonical shards is removed by~$\Theta$.
\end{prop}
\begin{proof}
The ``if'' direction follows from the definition of $\Psi/\Theta$ and the assertion of Proposition~\ref{recover cone} that $\Gamma$ is the intersection of its canonical shards.
The ``only if'' direction follows from Lemma~\ref{canon canon} and the observation that no shard of the form $\Sigma_\Theta(Y\covered Z)$, for $X\le Y\covered Z$, is removed by~$\Theta$, because $\Sigma_\Theta(Y\covered Z)$ separates two maximal cones of $\F/\Theta$.
\end{proof}

\begin{prop}\label{quotient shard face}
If $\Gamma\in\Psi/\Theta$ then any face of $\Gamma$ is in $\Psi/\Theta$.
\end{prop}
\begin{proof}
The cone $\Gamma\in\Psi/\Theta$ is in particular a cone in $\Psi$.
In the proof of Proposition~\ref{shard face}, we first showed that any facet~$C$ of any shard~$\Sigma$ is $\Gamma(F)$ for some codimension~2 face $F\in\F$.
Thus $C$ can be written as $\Sigma_1\cap\Sigma_2$, with each $\Sigma_i$ contained in the basic hyperplane $H_i$ of the rank-two full subarrangement consisting of hyperplanes containing~$C$.
Therefore $\Sigma_1\to\Sigma$ and $\Sigma_2\to\Sigma$ in the shard digraph, and thus by Proposition~\ref{good enough}, if~$\Theta$ removes $\Sigma_1$ or $\Sigma_2$ then~$\Theta$ must remove~$\Sigma$.
Thus any facet~$C$ of any unremoved shard $\Sigma\in\Psi/\Theta$ is an intersection of unremoved shards.

In the proof of Proposition~\ref{shard face}, we next showed that any facet~$F$ of a cone $\Gamma\in\Psi$ is the intersection of the canonical shards of $\Gamma$ with one additional shard $\Sigma'$.
This additional shard $\Sigma'$ defined a facet of one of the canonical shards of $\Gamma$, so as in the previous paragraph, if $\Sigma'$ is removed by~$\Theta$ then some canonical shard of $\Gamma$ is also removed.
Thus by Proposition~\ref{canon unremoved}, if $\Gamma\in\Psi/\Theta$ then $\Sigma'$ is not removed by~$\Theta$ so that~$F$ is the intersection of unremoved shards, i.e.\ $F\in\Psi/\Theta$.

The result for arbitrary faces follows just as in the proof of Proposition~\ref{shard face}.
\end{proof}

\begin{prop}\label{restrict bijection}
The map $\psi$ restricts to a bijection from $\pidown^\Theta(\R)$ to $\Psi/\Theta$.
The inverse map is the restriction of $\rho$ to $\Psi/\Theta$ and is an order-preserving map from $\Psi/\Theta$ to the restriction of the poset of regions to $\pidown^\Theta(\R)$.
\end{prop}
Recall from Proposition~\ref{bottoms} that the restriction of the poset of regions to $\pidown^\Theta(\R)$ is isomorphic to the lattice quotient $(\R,\le)/\Theta$.
\begin{proof}
In light of Proposition~\ref{bijection}, it is enough to show that $\psi$ maps $\pidown^\Theta(\R)$ into $\Psi/\Theta$ and that $\rho$ maps $\Psi/\Theta$ into $\pidown^\Theta(\R)$.
Suppose $R\in\pidown^\Theta(\R)$.
Then by Proposition~\ref{bottom char}, none of the canonical generators of~$R$ is contracted by~$\Theta$, or equivalently none of the shards associated to the canonical generators is removed.
Thus $\psi(R)$ is an intersection of unremoved shards.

On the other hand, let $\Gamma$ be an element of $\Psi/\Theta$, that is, an intersection of unremoved shards.
Then by Proposition~\ref{canon unremoved}, no canonical shard of $\Gamma$ is removed by~$\Theta$.
By Theorem~\ref{R can join rep}, none of the canonical joinands of $\rho(\Gamma)$ is contracted.
Thus by Proposition~\ref{bottom char}, $\rho(\Gamma)\in\pidown^\Theta(\R)$.
\end{proof}

By Proposition~\ref{restrict bijection} and the definition of the partial order $\preceq$, the map $\psi$ is an isomorphism between the restriction $(\pidown^\Theta(\R),\preceq)$ of $(\R,\preceq)$ and the poset $(\Psi/\Theta,\supseteq)$.
Since the join in $(\Psi,\supseteq)$ is intersection, it is immediate from the definitions that the induced subposet $(\Psi/\Theta,\supseteq)$ of $(\Psi,\supseteq)$ is a join-sublattice of $(\Psi,\supseteq)$.
Furthermore, $(\Psi/\Theta,\supseteq)$ is a lattice, since the bottom element $V$ of $(\Psi,\supseteq)$ is in $(\Psi/\Theta,\supseteq)$.
Thus 
\begin{prop}\label{cong sub}
The poset $(\pidown^\Theta(\R),\preceq)$ is a lattice and a join-sublattice of $(\R,\preceq)$.
\end{prop}

\begin{figure}
\scalebox{.8}{\includegraphics{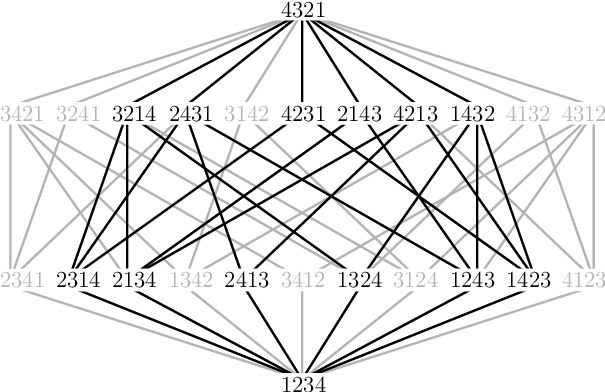}}
\caption{$(\pidown^{\Theta}(\R),\preceq)$ for the Cambrian congruence~$\Theta$ of Figure~\ref{A3camb_cong}.}
\label{S4 bip preceq}
\end{figure}

\begin{example}\label{NC in S4 preceq}
Figure~\ref{S4 bip preceq} shows the lattice $(\pidown^\Theta(\R),\preceq)$ in the case where $\A=\A(S_4)$ and~$\Theta$ is the congruence of Example~\ref{cong ex}.
See also Figure~\ref{S4 preceq}, Figure~\ref{A3camb_cong}, and Figure~\ref{Psi S4 bip}.
For easy comparison, $(\pidown^\Theta(\R),\preceq)$ is pictured superimposed on a faded view of the full lattice $(\R,\preceq)$.
\end{example}

We now generalize the properties of $(\R,\le)$ to $(\pidown^\Theta(\R),\le)$.

\begin{prop}\label{quotient graded}
The lattice $(\pidown^\Theta(\R),\preceq)$ is graded, with rank function equal to the number of lower hyperplanes of $R\in \pidown^\Theta(\R)$.
Alternately, the rank of a cone $\Gamma\in\Psi/\Theta$ is the codimension of $\Gamma$.
\end{prop}
\begin{proof}
Argue as in the proof of Proposition~\ref{graded}, with Proposition~\ref{Theta recover cone} replacing Proposition~\ref{recover cone}.
\end{proof}

\begin{prop}\label{quotient atomic coatomic}
The lattice $(\pidown^\Theta(\R),\preceq)$ is atomic and coatomic.
\end{prop}
\begin{proof}
Argue as in Proposition~\ref{atomic coatomic}, replacing Propositions~\ref{shard face} and~\ref{graded} by Propositions~\ref{quotient shard face} and~\ref{quotient graded}.
In the argument that $(\pidown^\Theta(\R),\preceq)$ is coatomic, we can assume that~$\Theta$ does not remove any of the shards which are basic hyperplanes in~$\A$.
If~$\Theta$ does remove a basic hyperplane $H$ of $\A$, then by Theorem~\ref{parabolic quotient}, we can realize $(\pidown^\Theta(\R),\preceq)$ as a subposet of $(\R(\A_{\B\setminus\set{H}}),\preceq)$.
\end{proof}

We do not generalize the first assertion of Proposition~\ref{psi fib}, since there seems to be no reasonable generalization of the intersection lattice associated to $\F/\Theta$.
However, the second assertion does generalize, by the identical proof, replacing Proposition~\ref{recover cone} by Proposition~\ref{Theta recover cone}.
\begin{prop}\label{quotient psi fib}
The lattice $(\Psi/\Theta,\supseteq)$ is anti-isomorphic to the fiber poset of $\Gamma_\Theta:(\F/\Theta)\to(\Psi/\Theta)$.
\end{prop}

The combinatorial description of $(\R,\preceq)$ generalizes as well.
For any $Q,R\in\pidown^\Theta(\R)$ with $Q\le R$, let 
\[I_\Theta(Q,R)=\set{J(P\covered P'):Q\le P\covered P'\le R,\,P'\in\pidown^\Theta(\R)}.\]
The order and cover relations $Q\le P\covered P'\le R$ refer to the poset of regions $(\R,\le)$.
Recall that $L(R)=\Meet\set{P:P\covered R}$.
The proof of Proposition~\ref{preceq combin}, with Proposition~\ref{Theta recover cone} replacing Proposition~\ref{recover cone}, proves the following:
\begin{prop}\label{quotient preceq combin}
If $Q,R\in\pidown^\Theta(\R)$ then $Q\preceq R$ if and only if $I_\Theta(\pidown^\Theta(L(Q)),Q)\subseteq I_\Theta(\pidown^\Theta(L(R)),R)$.
\end{prop}

The interval $[L(R),R]$ in $(\R(\A),\le_B)$ is isomorphic to $(\R(\A'),\le_{B'})$, where $\A'$ is the full subarrangement of~$\A$ consisting of hyperplanes $H\in\A$ with $H\supseteq\bigcap\Lower(R)$, and $B'$ is the $\A'$-region containing~$B$.
The interval $[L(R),R]$ is also the set of $\A$-regions containing $R\cap\bigcap\Lower(R)$, and the restriction of~$\Theta$ to $[L(R),R]$ is a lattice congruence on $[L(R),R]$.
This lattice congruence induces a congruence on $(\R(\A'),\le_{B'})$ which we will call $\Theta_R$.
The following is an immediate consequence of Propositions~\ref{lower int} and~\ref{restrict bijection}.

\begin{prop}\label{lower int quotient}
If $R\in\pidown^\Theta(\R)$, then the lower interval $[1,R]$ in \mbox{$(\pidown^\Theta(\R(\A)),\preceq_B)$} is isomorphic to $(\pidown^{\Theta_R}(\R(\A')),\preceq_{B'})$, where $\A'$, $B'$ and $\Theta_R$ are as above.
\end{prop}
In geometric terms, Proposition~\ref{lower int quotient} says that the lower interval below $\psi(R)$ in $(\Psi(\A,B)/\Theta,\supseteq)$ is isomorphic to $(\Psi(\A',B')/\Theta_R,\supseteq)$.

\begin{theorem}\label{quotient mobius}
Let~$\Theta$ be a congruence on $(\R,\le)$ and let $\B_\Theta$ be the set of basic hyperplanes of~$\A$ that are not removed by~$\Theta$.
Then the M\"{o}bius number of \mbox{$(\pidown^\Theta(\R),\preceq)$} is $(-1)^{|\B_\Theta|}$ times the number of regions in $\pidown^\Theta(\R)$ that are not contained in any proper parabolic subset of $\R_{\B_\Theta}$.
Equivalently, by inclusion-exclusion,
\[\mu_\preceq(B,\pidown^\Theta(-B))=\sum_{\K\subseteq \B_\Theta}(-1)^{|\K|}\left|\pidown^\Theta(\R_\K)\right|.\]
\end{theorem}

\begin{proof}
Let $\K\subseteq \B$ and let $F=B\cap\bigcap\K$.
By Proposition~\ref{facial}, the set of regions containing $F$ is an interval $[P,Q]$ in $(\R,\le)$.
Since~$B$ is in the interval, we have $P=B$.
By Lemma~\ref{RK}, $[B,Q]=\R_\K$.
Thus since $\pidown^\Theta(R)\le R$ for all~$R$, the quantity $\left|\pidown^\Theta(\R_\K)\right|$ is the number of congruence classes of the restriction of~$\Theta$ to $[B,Q]$.
We claim that the latter is in turn equal to the number of congruence classes of the restriction of~$\Theta$ to the interval $[-Q,-B]$ consisting of regions containing $-F=(-B)\cap\bigcap\K$.

To prove the claim, we first observe that $[B,Q]\cong[-Q,-B]$.
The isomorphism maps $P\in[B,Q]$ to the region whose separating set is $[\A\setminus S(Q)]\cup S(P)$.
The inverse maps $P\in[-Q,-B]$ to $S(Q)\setminus[\A\setminus S(P)]$.
Equivalently, the map and its inverse are $P\mapsto [(-Q)\join P]$ and $P\mapsto (P\meet Q)$.
The definition of congruence now allows us to easily conclude that the restriction of~$\Theta$ to $[B,Q]$ coincides with the restriction of~$\Theta$ to $[-Q,-B]$.
The claim follows.

For each $R\in\pidown^\Theta(\R)$, let $\A'(R)$ be as in the proof of Theorem~\ref{mobius}, let $R'$ be the $\A'(R)$-region containing~$R$, and let $\Theta_R$ be as in Proposition~\ref{lower int quotient}.
Then for any $\K\subseteq \Lower(R)$, the quantity $\left|\pidown^{\Theta_R}(\R(\A'(R))_\K)\right|$ is the number of congruence classes of the restriction of $\Theta_R$ to the interval $I'$ in $(\R(\A'(R)),\le)$ consisting of $\A'(R)$-regions containing the face $R'\cap\bigcap\K$.

As in the proof of Theorem~\ref{mobius}, let $G(R)$ be the intersection of all facets of~$R$ separating~$R$ from regions covered by~$R$ in $(\R,\le)$.
Recall that $\Theta_R$ is defined as the congruence on $(\R(\A'(R)),\le)$ corresponding the the restriction of~$\Theta$ to the interval of $\A$-regions containing $G(R)$. 
The interval $I$ consisting of $\A$-regions containing $R\cap\bigcap\K$ is a weakly smaller interval, so the number of classes of the restriction of~$\Theta$ to $I$ equals the number of classes of the restriction of $\Theta_R$ to $I'$.

Let $\Theta_{R,\K}$ be the restriction of~$\Theta$ to $I$.  
Then by Proposition~\ref{lower int quotient} and the previous paragraphs, the proposed M\"{o}bius function value on the interval $[B,R]$ in \mbox{$(\pidown^\Theta(W),\preceq)$} is 
\[\sum_{\K\subseteq \Lower(R)_{(\Theta_R)}}(-1)^{|\K|}\left(\#\mbox{ classes of }\Theta_{R,\K}\right),\]
where $\Lower(R)_{(\Theta_R)}$ is the set of lower hyperplanes of~$R$ (i.e.\ basic hyperplanes of $\A'(R)$) not removed by $\Theta_R$.
These correspond to lower shards of~$R$ not removed by~$\Theta$, but since $R\in\pidown^\Theta(\R)$, none of its lower shards is removed.
Thus the sum can be rewritten further as
\[\sum_{\K\subseteq \Lower(R)}(-1)^{|\K|}\left(\#\mbox{ classes of }\Theta_{R,\K}\right).\]

Let~$C$ be the maximal cone of $\F/\Theta$ containing~$R$ and let $G(C)$ be the intersection of all facets of~$C$ separating~$C$ from maximal cones covered by~$C$ in $(\R,\le)/\Theta$.
Recall from the end of Section~\ref{cong sec} that the lower hyperplanes of~$R$ are the hyperplanes which separate~$C$ from maximal cones covered by~$C$ in $(\R,\le)/\Theta$.
Thus $G(C)$ is the unique face in $\F/\Theta$ with $G(C)\supseteq G(R)$ and $\dim(G(C))=\dim(G(R))$.
Furthermore, the subsets $\K\subseteq\Lower(R)$ are in bijection with faces of~$C$ containing $G(C)$.

Recall also that $\F/\Theta$ is a coarsening of~$\F$, such that two $\A$-regions are in the same maximal cone of $\F/\Theta$ if and only if they are congruent modulo~$\Theta$.
Thus the classes in the restriction $\Theta_{R,\K}$ of~$\Theta$ to regions containing $R\cap\bigcap\K$ are in bijection with the maximal cones in $\F/\Theta$ contining $C\cap\bigcap\K$.
Now the proposed M\"{o}bius number for $[B,R]$ is 
\[\sum_{\substack{F\in\F/\Theta\\C\supseteq F\supseteq G(C)}}(-1)^{\codim(F)}\left(\#\mbox{ maximal faces }C'\in\F/\Theta\mbox{ with }C'\supseteq F\right).\]
We replace the sum over $\pidown^\Theta(\R)$ by a sum over maximal cones $C$ of $\F/\Theta$ and complete the proof just as the proof of Theorem~\ref{mobius}, replacing Theorem~\ref{shelling} by the analogous fact for $\F/\Theta$.
\end{proof}

Let~$\Theta$ be some congruence on $(\R,\le)$ and let $\S(\F/\Theta)$ be the dual sphere to $\F/\Theta$.
Let $\overline{\S}(\F/\Theta)$ be the regular CW-ball obtained by gluing onto $\S(\F/\Theta)$ a cell of dimension one higher than the dimension of $\S(\F/\Theta)$.
Let $\Delta(\F/\Theta)$ be the pulling triangulation of $\overline{\S}(\F/\Theta)$ obtained by ordering the vertices of $\overline{\S}(\F/\Theta)$ according to the lattice quotient $(\R,\le)/\Theta$.
In the case where $\F/\Theta$ is the normal fan of some polytope $P$, $\Delta(\F/\Theta)$ is a triangulation of $P$.
Now Theorem~\ref{zon chain tri} and Proposition~\ref{restrict bijection} imply the following theorem.

\begin{theorem}\label{chain tri cong}
The map $\delta$ restricts to a dimension-preserving bijection between the order complex of $(\pidown^\Theta(\R),\preceq)$ and the triangulation $\Delta(\F/\Theta)$, with inverse map~$\gamma$.
\end{theorem}

\begin{figure}
\includegraphics{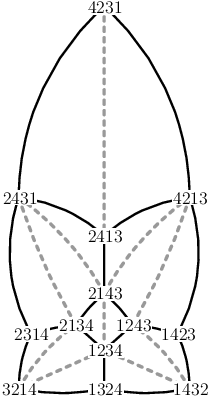}
\caption{A stage in the construction of $\Delta(\F/\Theta)$.}
\label{A3 bip pull fig}
\end{figure}

\begin{example}\label{chain tri cong example}
Figure~\ref{A3 bip pull fig} depicts a stage in the construction of $\Delta(\F/\Theta)$ in the case where $\A=\A(S_4)$ and~$\Theta$ is the congruence pictured in Figure~\ref{A3camb_cong}, with drawing conventions similar to Figure~\ref{A3 pull fig}.
Theorem~\ref{chain tri cong} says that the 16 maximal chains of $(\Psi/\Theta,\preceq)$ correspond to (cones over the) 16 triangles shown in the figure.
\end{example}

\begin{remark}\label{subcomplex}
For the congruence~$\Theta$ of Example~\ref{chain tri cong example}, we see that $\Delta(\F/\Theta)$ is an induced subcomplex of $\Delta(\A)$.
This property does not hold in general.
For example, there is a lattice congruence $\Theta'$ on the weak order on $S_4$ that contracts the edge $2143\covered2413$ and no other edges.
Aided by Figure~\ref{A3 pull fig}, it is easy to see that $\Delta(\F/\Theta')$ is not a subcomplex of $\Delta(\A)$.
We do not know how general this subcomplex property is, but it is natural to wonder whether it holds for all of the Cambrian congruences, defined in the next section.
\end{remark}

\section{Noncrossing partition lattices}\label{nc sec}
In this section, we consider the lattices $(\pidown^\Theta(W),\preceq)$ in the special case where~$\Theta$ is a Cambrian congruence.

A \emph{Coxeter element}~$c$ of~$W$ is the product, in any order, of the simple generators $S$ of~$W$.
We fix some particular Coxeter element $c=s_1s_2\cdots s_n$, for $S=\set{s_1,\ldots,s_n}$.
Let $m$ be the function such that the defining relations of~$W$ include $(s_is_j)^{m(s_i,s_j)}=1$ for all $s_i\neq s_j\in S$.
We define the \emph{Cambrian congruence}~$\Theta_c$ by specifying a set of join-irreducible elements that~$\Theta_c$ is required to contract.
For every pair $s_i,s_j\in S$ with $i<j$, we require that~$\Theta_c$ contracts all join-irreducible elements having a reduced word alternating $s_js_is_js_i\cdots$ of length at least two and at most $m(s_i,s_j)-1$.
The $c$-Cambrian congruence~$\Theta_c$ is the finest congruence contracting those join-irreducible elements.
The \emph{$c$-Cambrian lattice} is the quotient of the weak order on~$W$ modulo the $c$-Cambrian congruence.

\begin{example}\label{H4 camb}
Let~$W$ be the Coxeter group (of type~$H_4$) with $S=\set{q,r,s,t}$, $m(q,r)=5$, $m(r,s)=m(s,t)=3$ and $m(q,s)=m(q,t)=m(s,t)=2$.  Choose the Coxeter element $c=rqts$.
The Cambrian congruence~$\Theta_c$ is then the finest lattice congruence that contracts $qr$, $qrq$, $qrqr$, $sr$ and $st$.
\end{example}

\begin{example}\label{S4 camb}
When~$W$ is $S_4$ with $S=\set{(1\,\,2),(2\,\,3),(3\,\,4)}$ and $c=(1\,\,2)(3\,\,4)(2\,\,3)$, the Cambrian congruence~$\Theta_c$ is the congruence pictured in Figure~\ref{A3camb_cong}.
That is, the congruence pictured in Figure~\ref{A3camb_cong} is the finest congruence on the weak order on $S_4$ contracting $(2\,\,3)(1\,\,2)=3124$ and $(2\,\,3)(3\,\,4)=1342$, or in other words, the finest congruence with $1324\equiv3124$ and $1324\equiv1342$.

We relate this example to several previous examples.
One verifies in Figure~\ref{A3camb_cong} that~$\Theta_c$ contracts the join-irreducible elements $3124$, $1342$, $2341$, $4123$, $3412$ and no other join-irreducible elements.
This is in keeping with Figure~\ref{A3digraph fig}, where we see that the smallest order ideal containing $\set{\Sigma(3124),\Sigma(1342)}$ in the transitive closure of the shard digraph is indeed $\set{\Sigma(3124),\Sigma(1342),\Sigma(2341),\Sigma(4123),\Sigma(3412)}$.
\end{example}

For each $c$, write $\pidown^c$ for $\pidown^{\Theta_c}$ and $\piup_c$ for $\piup_{\Theta_c}$.
The elements of the set $\pidown^c(W)$ are called \emph{$c$-sortable elements}.
As a special case of Proposition~\ref{bottoms}, the $c$-Cambrian lattice is (isomorphic to) the weak order on $c$-sortable elements.
In~\cite{sortable}, $c$-sortable elements are defined in terms of the combinatorics of reduced words, and in~\cite{sort_camb}, the $c$-sortable elements are shown to be the bottom elements of congruence classes of~$\Theta_c$.
Here, we do not need the definition of $c$-sortable elements in terms of reduced words.

\begin{example}\label{S4 sort}
When~$W$ is $S_4$ and $c=(1\,\,2)(3\,\,4)(2\,\,3)$, the $c$-sortable elements of~$W$ are the permutations appearing in Figures~\ref{Psi S4 bip} and~\ref{S4 bip preceq}.
As a special case of a characterization of the inversion sets of $c$-sortable elements (see \cite[Theorem~4.1]{sortable} or \cite[Theorem~4.3]{typefree}), these are the permutations in $S_4$ not containing $312$, $412$, $342$ or $341$ as a subsequence.
\end{example}

Let~$T$ be the set of reflections of~$W$.
A \emph{reduced~$T$-word} for $w\in W$ is a shortest possible word for~$w$ in the alphabet~$T$.
(This contrasts with the usual notion of a \emph{reduced word} for~$W,$ a shortest possible word for~$w$ in the alphabet $S$.)
The \emph{absolute order} on~$W$ is the prefix order on reduced~$T$-words:  we set $u\le v$ if and only if every reduced~$T$-word for $u$ occurs as a prefix of some reduced~$T$-word for $v$.
The \emph{$W$-noncrossing partition lattice} $\NC_c(W)$ is the interval $[1,c]_T$ in the absolute order, where~$c$ is a Coxeter element of~$W\!$.
It is straightforward to show that the isomorphism type of $\NC_c$ does not depend on the choice of~$c$.
The elements of $[1,c]_T$ are called \emph{$c$-noncrossing partitions}.

For the present purposes, the most important result about $c$-sortable elements and $c$-noncrossing partitons is the following theorem, which is \cite[Theorem~6.1]{sortable}.  
(Cf. \cite[Theorem~8.9]{typefree}.)
Recall that the descents of~$w$ are the simple generators $s\in S$ such that $\ell(ws)<\ell(w)$.
\begin{theorem}\label{nc}
For any Coxeter element~$c$, there is a bijection $w\mapsto\nc_c(w)$ from the set of $c$-sortable elements to the set of $c$-noncrossing partitions.
Furthermore $\nc_c$ maps $c$-sortable elements with $k$ descents to $c$-noncrossing partitions of rank $k$.
\end{theorem}

The lattice $(\pidown^c(W),\preceq)$ is the restriction of $(W,\preceq)$ to the $c$-sortable elements of~$W$.
The main results of this section are the following theorem and corollary.

\begin{theorem}\label{isom}
The map $\nc_c$ is an isomorphism from $(\pidown^c(W),\preceq)$ to $\NC_c(W)$.
\end{theorem}

\begin{cor}\label{lattice}
The poset $\NC_c(W)$ is a lattice.
\end{cor}

Corollary~\ref{lattice} follows immediately from Theorem~\ref{isom} and Proposition~\ref{cong sub}, which states that $(\pidown^\Theta(W),\preceq)$ is a lattice and a join-sublattice of $(W,\preceq)$ for any congruence~$\Theta$.
In fact, more is true in the case $\Theta=\Theta_c$.
Recall that given a cone $\Gamma\in\Psi$, the notation $U(\Gamma)$ represents the intersection of all hyperplanes in~$\A$ containing $\Gamma$.

\begin{prop}\label{sublattice}
The lattice $(\pidown^c(W),\preceq)$ is a sublattice of $(W,\preceq)$.
\end{prop}

\begin{prop}\label{meet sublattice}
The map~$U$ embeds $(\Psi/\Theta_c,\supseteq)$ as a meet-sublattice of $\Int(\A)$.
\end{prop}

Before proving Theorem~\ref{isom} and Propositions~\ref{sublattice} and~\ref{meet sublattice}, we discuss the consequences that follow from Theorem~\ref{isom} and the results of Section~\ref{cong psi sec}.
Corollary~\ref{nc mobius} is the concatenation of \cite[Corollary 7.4.ii]{ABMW} and \cite[Corollary~4.4]{ABW}.

\begin{theorem}\label{Thetac mobius}
The M\"{o}bius function of $(\pidown^{\Theta_c}(W),\preceq)$ is $(-1)^{\rank(W)}$ times the number of elements of $\pidown^c(W)$ that are not contained in any proper standard parabolic subgroup of~$W$.
\end{theorem}
\begin{cor}\label{nc mobius}
The M\"{o}bius number $\mu_{\le_T}(1,c)$ of $\NC_c(W)$ is $(-1)^{\rank(W)}$ times the number of elements of $[1,c]_T$ not contained in any proper standard parabolic subgroup of~$W$.
\end{cor}

\begin{prop}\label{nc inv op}
The inverse of $\nc_c$ is an order-preserving map from $\NC_c(W)$ to the $c$-Cambrian lattice.
\end{prop}
Proposition~\ref{nc inv op} says that $\NC_c(W)$ is a weaker partial order than the $c$-Cambrian lattice.
This fact has never been published, but has been observed independently by several researchers in the brief time since the results of~\cite{cambrian} were disseminated.

When~$\Theta$ is a Cambrian congruence~$\Theta_c$, the fan $\F/\Theta_c$ is the normal fan of a $W$-associahedron \cite[Theorem~3.4]{HoLaTh}, so $\overline{\S}(\F/\Theta_c)$ is a $W$-associahedron.
Extending $\nc_c$ in the natural way to a bijection on chains, we have the following immediate consequence of Theorem~\ref{chain tri cong}.

\begin{theorem}\label{chain tri camb}
The map $\delta\circ\nc_c^{-1}$ is a dimension-preserving bijection between the order complex of the noncrossing partition lattice $\NC_c(W)=[1,c]_T$ and the pulling triangulation $\Delta(\F/\Theta_c)$ of the $W$-associahedron.
\end{theorem}
The pulling triangulation $\Delta(\F/\Theta_c)$ is obtained by pulling the vertices of the $W$-associahedron in the dual order to the Cambrian lattice $\pidown^{\Theta_c}(W)$.

\begin{example}\label{camb tri example}
In Example~\ref{chain tri cong example} and Figure~\ref{A3 bip pull fig}, the congruence~$\Theta$ is the Cambrian congruence of Example~\ref{S4 camb}.
\end{example}

\begin{remark}\label{Loday remark}
We explain the relationship between Theorems~\ref{zon chain tri}, \ref{chain tri cong}, and~\ref{chain tri camb} and the construction given in~\cite{Loday}.
In the latter reference, Loday constructed the pulling triangulation of the classical associahedron with vertices ordered by the Tamari lattice.
He gave a bijection between parking functions and maximal simplices of this pulling triangulation.
He also considered the weak-order pulling triangulation of the classical permutohedron and asked whether the maximal simplices are in bijection with some object analogous to parking functions.

Since parking functions are in bijection with maximal chains in the usual noncrossing partition lattice $\NC_c(S_n)$, one might expect, for any~$W$, the maximal simplices in the pulling triangulation of the $W$-associahedron to be in bijection with maximal chains in $\NC_c(W)$.
Theorem~\ref{chain tri camb} says that much more is true.
Furthermore, Theorem~\ref{zon chain tri} shows that the analogous relationship holds between $(W,\preceq)$ and the pulling triangulation of the $W$-permutohedron, and, more generally, between $(\R,\preceq)$ and the zonotope associated to $\A$.
Theorem~\ref{chain tri cong} generalizes Theorem~\ref{zon chain tri} to a broad level of generality that includes Theorem~\ref{chain tri camb}.
\end{remark}

We now prepare to prove Theorem~\ref{isom}. 
Descents of $w$ are in bijection with \emph{cover reflections} of~$w$:
For each descent~$s$ of $w$, the reflection $wsw^{-1}$ is a cover reflection.
Let $\cov(w)$ denote the cover reflections of $w$.
Recall that~$H_t$ is the reflecting hyperplane for the reflection~$t$.
The lower hyperplanes of the $\A(W)$-region associated to $w$ are $\set{H_t:t\in\cov(w)}$.

An immediate corollary \cite[Corollary~3.9]{camb_fan} to Theorem~\ref{nc} is that for each reflection~$t$ of~$W$, there exists a unique $c$-sortable join-irreducible element of the weak order whose unique cover reflection is~$t$.
Combining this with Proposition~\ref{shard ji}, we have the following proposition:

\begin{prop}\label{one}
Each hyperplane of $\A(W)$ contains exactly one shard that is not removed by~$\Theta_c$.
\end{prop}

We write $\F_c$ for the fan $\F/\Theta_c$ and call $\F_c$ the \emph{$c$-Cambrian fan}.
For any cone $F\in V$ and any $x\in F$, the \emph{linearization $\Lin_x(F)$ of~$F$ at $x$} is the set of vectors $x'\in V$ such that $x+\epsilon x'$ is in~$F$ for any sufficiently small nonnegative~$\epsilon$.
If~$\F$ is a fan and $x$ is any point in the relative interior of some cone of~$\F$ then the \emph{star} of $x$ in~$\F$ is the fan $\Star_x(\F)=\set{\Lin_x(F): x \in F\in\F}$.  
For any cone~$F$ of~$\F$, the \emph{star} of~$F$ in~$\F$ is the fan $\Star_F(\F)=\Star_x(\F)$ for any $x$ in the relative interior of~$F$.
The following is part of \cite[Proposition~8.12]{camb_fan}.

\begin{prop}\label{star camb}
For any face $F\in\F_c$, the star of~$F$ in $\F_c$ is isomorphic to a Cambrian fan associated to the parabolic subgroup $W'$ generated by reflections fixing~$F$ pointwise.
\end{prop}

Recall that we fixed a representation of~$W$ as a group of orthogonal transformations of the Euclidean vector space~$V$.
Let $\Fix:W\to\Int(\A(W))$ be the map taking $w\in W$ to the subspace of~$V$ consisting of points fixed by $w$.
The following facts are proved in~\cite{BWorth} and \cite[Section~2]{BWKpi}:
The restriction of the map $\Fix$ to $[1,c]_T$ is injective; we call the subspaces in the image $\Fix([1,c]_T)$ the \emph{$c$-noncrossing subspaces}.
The map $\Fix$ is an isomorphism from $[1,c]_T$ to the poset of $c$-noncrossing subspaces under reverse containment.

\begin{proof}[Proof of Theorem~\ref{isom}]
The map $\Fix\circ\nc_c$ takes $w\in\pidown^c(W)$ to the intersection of the hyperplanes associated to $\cov(w)$.
This map coincides with the restriction of $U\circ\psi$, where $\psi$ is the map from Proposition~\ref{bijection} and~$U$ is the map from Proposition~\ref{psi fib}.
Now $\nc_c$ is a bijection from $\pidown^c(W)$ to $[1,c]_T$, the map $\Fix$ restricts to an isomorphism between $[1,c]_T$ and the $c$-noncrossing subspaces under reverse containment, and $\psi$ restricts to an isomorphism from $(\pidown^c(W),\preceq)$ to $(\Psi/\Theta_c,\supseteq)$.
Thus~$U$ restricts to a bijection from $\Psi/\Theta_c$ to the $c$-noncrossing subspaces.
This restriction is containment-preserving, and we can complete the proof by showing that its inverse is also containment-preserving.

Suppose $U_1$ and $U_2$ are noncrossing subspaces with $U_1\supseteq U_2$.
Let~$\Gamma_1$ and~$\Gamma_2$ be the preimages of $U_1$ and $U_2$ in $\Psi/\Theta_c$.
To show that $\Gamma_1\supseteq \Gamma_2$, we show that every shard containing~$\Gamma_1$ also contains~$\Gamma_2$.
Let~$\Sigma$ be a shard containing~$\Gamma_1$.
By Proposition~\ref{one},~$\Sigma$ is the only shard contained in $H(\Sigma)$ that is not removed by~$\Theta_c$.
Since the map~$U$ is order-preserving and $\Sigma\supseteq \Gamma_1\supseteq \Gamma_2$, we have $H(\Sigma)=U(\Sigma)\supseteq U(\Gamma_2)=U_2$.
Let~$F$ be a face in $\F_c$ such that~$F\subseteq\Gamma_2$ and $\dim(F)=\dim(\Gamma_2)$.
By Proposition~\ref{star camb}, the star of~$F$ in $\F_c$ is isomorphic to a Cambrian fan associated to a parabolic subgroup of~$W$.
Thus Proposition~\ref{one} implies that, for each hyperplane~$H$ in~$\A$ containing~$F$, there is a shard in~$H$ containing~$F$ and not removed by~$\Theta_c$.
But $H(\Sigma)\supseteq U_2\supseteq F_2$ and~$\Sigma$ is the only unremoved shard in $H(\Sigma)$, so $\Sigma\supseteq F$.
By Proposition~\ref{recover cone},~$\Gamma_2$ is the intersection of all shards containing~$F$, so $\Sigma\supseteq \Gamma_2$.
\end{proof}

\begin{proof}[Proof of Proposition~\ref{meet sublattice}]
Suppose $U_1$ and $U_2$ are $c$-noncrossing subspaces and let~$U_3$ be their meet in $\Int(\A)$ (not \emph{a priori} the meet in the lattice of $c$-noncrossing subspaces).
Thus~$U_3$ is the intersection of all hyperplanes in~$\A$ containing both $U_1$ and $U_2$.
Let~$\Gamma_1$ and~$\Gamma_2$ be the preimages of $U_1$ and $U_2$ in $\Psi/\Theta_c$.
The meet of~$\Gamma_1$ and~$\Gamma_2$ in $(\Psi/\Theta_c,\supseteq)$ is the intersection of all unremoved shards containing both~$\Gamma_1$ and~$\Gamma_2$.
Thus to prove the proposition, it is enough to show that for each hyperplane~$H$ containing $U_1$ and $U_2$, the unique unremoved shard~$\Sigma$ contained in~$H$ has $\Sigma\supseteq \Gamma_1$ and $\Sigma\supseteq \Gamma_2$.
Arguing exactly as in the proof of Theorem~\ref{isom}, we conclude that $\Sigma\supseteq \Gamma_2$, and repeating the argument, we conclude that $\Sigma\supseteq \Gamma_1$.
\end{proof}

\begin{proof}[Proof of Proposition~\ref{sublattice}]
Let~$\Gamma_1$ and~$\Gamma_2$ be cones in $\Psi/\Theta_c$.
Then by Proposition~\ref{meet sublattice}, the rank of the meet of~$\Gamma_1$ and~$\Gamma_2$ in $(\Psi/\Theta_c,\supseteq)$ equals the rank of the meet of $U(\Gamma_1)$ and $U(\Gamma_2)$ in $\Int(\A)$.
But since~$U$ is a rank-preserving and order-preserving map from $(\Psi,\supseteq)$ to $(\Int(\A),\supseteq)$, the rank of the meet of $U(\Gamma_1)$ and $U(\Gamma_2)$ in $\Int(\A)$ is weakly greater than the rank of the meet of~$\Gamma_1$ and~$\Gamma_2$ in $(\Psi,\supseteq)$.
Thus, since $(\Psi/\Theta_c,\supseteq)$ is an induced subposet of $(\Psi,\supseteq)$, the meet of~$\Gamma_1$ and~$\Gamma_2$ must be the same in $(\Psi/\Theta_c,\supseteq)$ as in $(\Psi,\supseteq)$.
\end{proof}

\begin{remark}\label{BW remark}
We now compare the proof of the lattice property (Corollary~\ref{lattice}) to the previously known proof from~\cite{BWlattice}.
A generator in $S$ is \emph{initial} in a Coxeter element~$c$ if there is some reduced word for~$c$ whose first letter is~$s$.
Similarly,~$s$ is \emph{final} in~$c$ if~$s$ is the last letter of some reduced word for~$c$.
Let $C$ be the cone consisting of points weakly above all hyperplanes~$H_s$ for~$s$ initial in~$c$ and weakly below all hyperplanes~$H_s$ for~$s$ final in~$c$. 
Let $\widetilde{\Psi}/\Theta_c=\set{\Gamma\cap C:\Gamma\in\Psi/\Theta_c}$.
Using results of \cite{camb_fan,typefree}, one can show that $(\widetilde{\Psi}/\Theta_c,\supseteq)$ is isomorphic to $(\Psi/\Theta_c,\supseteq)$.
For any~$W$, it is possible to choose a \emph{bipartite} Coxeter element $c$, meaning that every $s\in S$ is either initial or final in~$c$.
In this case $\widetilde{\Psi}/\Theta_c$ is contained in a simplicial cone whose facets are defined by the hyperplanes~$H_s$ for $s\in S$.

\begin{example}\label{BW example}
The Coxeter element $c=(1\,\,2)(3\,\,4)(2\,\,3)$ in $S_4$ is bipartite.
The collection $\widetilde{\Psi}/\Theta_c$ is shown in Figure~\ref{Psi A3 bip part}.
The cones illustrated in the figure are the intersections of the cones in Figure~\ref{Psi A3} with the closed cone consisting of all points (weakly) above $H_{(1\,2)}$ and $H_{(3\,4)}$ and below $H_{(2\,3)}$.
\end{example}

\begin{figure}
\scalebox{.5}{\includegraphics{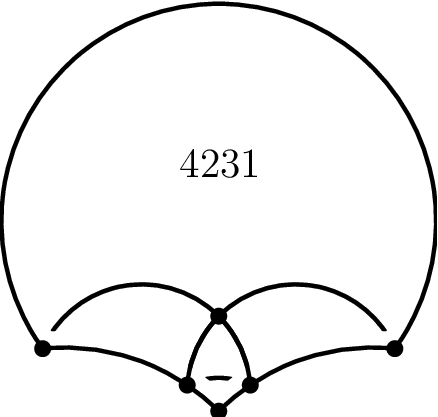}}
\caption{$\widetilde{\Psi}/\Theta_c$ for $W=S_4$ and $c=(1\,\,2)(3\,\,4)(2\,\,3)$.}
\label{Psi A3 bip part}
\end{figure}

The collection $\widetilde{\Psi}/\Theta_c$ is closely related to the complex $X(\gamma)$ defined in~\cite{BWlattice}.
(The symbol $\gamma$ denotes the Coxeter element called~$c$ here.)
The proof in~\cite{BWlattice} that $\NC_c(W)$ is a lattice proceeds as follows: one identifies coatoms of $\NC_c(W)$ with certain codimension-1 cones in $X(\gamma)$ and then shows that $\NC_c(W)$ is the poset of intersections these codimension-1 cones.
In particular, meets exist in $\NC_c(W)$ as intersections, so that $\NC_c(W)$ is a lattice.
In the present paper, we have identified the atoms of $\NC_c(W)$ with (codimension-one) cones (pieces $\widetilde{\Sigma}$ of shards), and shown that $\NC_c(W)$ is the poset of intersections of shards under reverse containment.
In particular, joins exist in $\NC_c(W)$ as intersections, so that $\NC_c(W)$ is a lattice.
There is a complete fan $AX(\gamma)$ (see \cite[Theorem~8.3]{BWlattice}) related to $X(\gamma)$ as the Cambrian fan $\F_c$ is related to the collection of shards unremoved by~$\Theta_c$.
The two fans are linearly isomorphic. 
(See \cite[Theorem~9.1]{camb_fan} or \cite[Theorem~5.5]{BWpermassoc}.  The paper~\cite{BWpermassoc} explains in detail the relationship between $AX(\gamma)$ and the Cambrian fan.)
This linear isomorphism, together with the self-duality of $\NC_c(W)$, suggests a duality between the broadest outlines (but not the details) of the two approaches to proving the lattice property. 
\end{remark}

\begin{remark}\label{semiinvariant remark}
Choosing a Coxeter element $c$ of a finite Coxeter group $W$ is equivalent to choosing an orientation of the Coxeter diagram of $W$, and this orientation is a \emph{Dynkin quiver} $Q$ when $W$ is of type $A$, $D$, or $E$.
There is a connection between the set $\Psi/\Theta_c$ and the \emph{domains of virtual semi-invariants} considered in~\cite{IOTW}.
The connection relies on known results which would require a substantial quiver-theoretic digression to explain.
Here we merely mention the connection, with the hope of making an exposition of the details in a future paper.

For a Dynkin quiver, the maximal domains of virtual semi-invariants are certain codimension-1 cones.
There is an invertible linear map which takes the set of maximal domains of virtual semi-invariants to the set of shards not removed by $\Theta_c$.
One can show further that this linear map induces an isomorphism between $(\Psi/\Theta_c,\supseteq)$ and the reverse inclusion order on domains of virtual semi-invariants.
Thus the reverse inclusion order on domains of virtual semi-invariants is isomorphic to $\NC_c(W)$.
For finite (crystallographic) Coxeter groups of other types, a similar statement holds via a standard \emph{folding} argument.

In light of such a fundamental quiver-theoretic interpretation of the shards not removed by $\Theta_c$, it is natural to ask whether there is a quiver-theoretic interpretation of the full set of shards in the Coxeter arrangement.
Furthermore, results of~\cite{Chindris} extend the consideration of virtual semi-invariants beyond quivers of finite type, raising the possibility of defining some generalization of ``shards not removed by $\Theta_c$'' to infinite Coxeter groups.
The generalization of unremoved shards would reproduce, and perhaps extend, some of the results obtained in~\cite{typefree} on Cambrian fans for infinite Coxeter groups.
\end{remark}

\begin{remark}\label{lattice remark}
The \emph{degree} of a join-irreducible region $J$ is the smallest positive integer $d$ such that $J$ is contained in some standard parabolic subset $\R_\K$ with $|\K|=d$.
A congruence~$\Theta$ has degree at most $d$ if it is generated by contracting a collection of join-irreducible regions, each of which has degree at most $d$.
Using Lemmas~\ref{whole} and~\ref{half} and Proposition~\ref{shard face}, one can fairly easily prove the following proposition, which shows in particular to what extent Proposition~\ref{sublattice} is special.

\begin{prop}\label{degree 2}
Suppose the shard digraph associated to $(\A,B)$ is acyclic.
Let~$\Theta$ be a lattice congruence on $(\R,\le)$.
If $(\pidown^\Theta(\R),\preceq)$ is a sublattice of $(\R,\preceq)$ then~$\Theta$ is of degree at most two.
\end{prop}
\end{remark}

\section*{Acknowledgments}
The author thanks David Speyer and Calin Chindris for bringing to his attention the connection between semi-invariants and Cambrian fans, and Marcelo Aguiar and Christos Athanasiadis for helpful comments on an earlier version of this paper.
This research was greatly assisted by Stembridge's \texttt{maple} packages \texttt{coxeter/weyl} and \texttt{posets}, and by Sloane's On-Line Encyclopedia of Integer Sequences.
An extended abstract of this paper was submitted to the 2009 FPSAC conference, and thus appears as~\cite{FPSAC}.  
Some parts of the extended abstract are direct quotes from this paper.
The author was partially supported by NSA grant H98230-09-1-0056.


\begin{thebibliography}{99}

\bibitem{AgSo}
M. Aguiar and F. Sottile,
\textit{Structure of the Malvenuto-Reutenauer Hopf algebra of permutations.}
Adv. Math. \textbf{191} (2005), no. 2, 225--275. 

\bibitem{ABMW}
C.~A.~Athanasiadis, T.~Brady, J.~McCammond, and C.~Watt,
\textit{$h$-vectors of generalized associahedra and non-crossing partitions.}
Int.\ Math.\ Res.\ Not.\ \textbf{2006}, Art.~ID 69705, 28 pp. 

\bibitem{ABW}
C. A. Athanasiadis, T. Brady and C. Watt,
\textit{Shellability of noncrossing partition lattices.}
Proc. Amer. Math. Soc. \textbf{135} (2007), no. 4, 939--949. 

\bibitem{Ath-Rei}
C. A. Athanasiadis and V. Reiner,
\textit{Noncrossing partitions for the group $D_n$.}
SIAM J. Discrete Math. \textbf{18} (2004), no.~2, 397--417.

\bibitem{Bessis}
D. Bessis,
\textit{The dual braid monoid.}
Ann. Sci. \'{E}cole Norm. Sup. (4) \textbf{36} (2003), no.~5, 647--683.

\bibitem{Biane1}
P. Biane,
\textit{Some properties of crossings and partitions.}
Discrete Math. \textbf{175} (1997), no.~1-3, 41--53.


\bibitem{BEZ}
A. Bj\"{o}rner, P. Edelman and G. Ziegler,
\textit{Hyperplane Arrangements with a Lattice of Regions.}
Discrete Comput. Geom. \textbf{5} (1990), 263--288.


\bibitem{BWorth}
T. Brady and C. Watt,
\textit{A partial order on the orthogonal group.}
Comm. Algebra \textbf{30} (2002) no.~8, 3749--3754.

\bibitem{BWKpi}
T. Brady and C. Watt,
\textit{$K(\pi,1)$'s for Artin groups of finite type.}
Proceedings of the Conference on Geometric and Combinatorial Group Theory, Part I (Haifa, 2000).
Geom. Dedicata \textbf{94} (2002), 225--250.

\bibitem{BWlattice}
T.~Brady and C.~Watt, 
\textit{Non-crossing partition lattices in finite real reflection groups.}
Trans.\ Amer.\ Math.\ Soc., \textbf{360} (2008), 1983--2005. 

\bibitem{BWpermassoc}
T.~Brady and C.~Watt, 
\textit{From Permutahedron to Associahedron.}
Preprint, 2008 (arXiv:0804.2331).

\bibitem{Cha-Sn}
I. Chajda and V. Sn\'a\v{s}el,
\textit{Congruences in Ordered Sets.}
Math. Bohem. \textbf{123} (1998) no. 1, 95--100.

\bibitem{Chapoton}
F. Chapoton,
\textit{Enumerative properties of generalized associahedra.}
S\'{e}m. Lothar. Combin. \textbf{51} (2004/05), Art. B51b, 16 pp.

\bibitem{Chindris}
C. Chindris,
\textit{Cluster fans, stability conditions, and domains of semi-invariants.}
Trans.\ Amer.\ Math.\ Soc., to appear.

\bibitem{gcccc}
S. Fomin and N. Reading,
\textit{Generalized cluster complexes and Coxeter combinatorics.}
Int. Math. Res. Not. \textbf{2005}, no.~44, 2709--2757. 

\bibitem{ga}
S. Fomin and A. Zelevinsky,
\textit{$Y$-systems and generalized associahedra.}
Ann. of Math. (2) \textbf{158} (2003), no. 3, 977--1018.

\bibitem{FreeLattices}
R. Freese, J. Je\v{z}ek, and J. Nation,
\textit{Free lattices.} 
Mathematical Surveys and Monographs \textbf{42}. 
American Mathematical Society, Providence, RI, 1995.

\bibitem{Gratzer}
G. Gr\"{a}tzer,
General Lattice Theory.
Second edition,
Birkh\"{a}user Verlag, Basel, 1998.

\bibitem{HoLaTh}
C. Hohlweg, C. Lange and H. Thomas,
\textit{Permutohedra and generalized associahedra.}
Preprint, 2007 (arXiv:0709.4241).

\bibitem{IOTW}
K. Igusa, K. Orr, G. Todorov and J. Weyman.
\textit{Cluster Complexes via Semi-Invariants.}
Preprint, 2007 (arXiv:0708.0798).

\bibitem{IngThom}
C. Ingalls and H. Thomas,
\textit{Noncrossing partitions and representations of quivers.}
preprint, 2006 (\texttt{math.RT/0612219}).

\bibitem{Jed}
P. Jedli\v{c}ka,
\textit{A Combinatorial Construction of the Weak Order of a Coxeter Group.}
Comm. in Alg. \textbf{33} (2005), 1447--1460.

\bibitem{Kreweras}
G. Kreweras,
\textit{Sur les partitions non crois\'{e}es d'un cycle.}
Discrete Math. \textbf{1} (1972), no.~4, 333--350.

\bibitem{Lee}
C. W. Lee,
\textit{Subdivisions and triangulations of polytopes}. 
Handbook of discrete and computational geometry, 271--290, 
CRC Press Ser. Discrete Math. Appl., CRC, Boca Raton, FL, 1997. 

\bibitem{Loday}
J.-L. Loday,
\textit{Parking functions and triangulation of the associahedron.} 
Categories in algebra, geometry and mathematical physics, 327--340, 
Contemp.\ Math. \textbf{431}, Amer. Math. Soc., Providence, RI, 2007. 

\bibitem{Postnikov}
A.~Postnikov,
\textit{Permutohedra, associahedra and beyond.}
preprint, 2005 (\texttt{arXiv:math/0507163}).

\bibitem{recursion}
N.~Reading, 
\textit{The cd-index of Bruhat intervals.} 
Electron.\ J.\ Combin.\ \textbf{11(1)} (2004), Research Paper 74, 25 pp. (electronic).

\bibitem{hyperplane}
N. Reading,
\emph{Lattice and order properties of the poset of regions in a hyperplane arrangement.}
Algebra Universalis \textbf{50} (2003), 179--205.

\bibitem{hplanedim}
N.~Reading,
\textit{The order dimension of the poset of regions in a hyperplane arrangement.}
J.\ Combin.\ Theory Ser.\ A, \textbf{104} (2003) no.~2, 265--285.

\bibitem{congruence}
N.~Reading,
\textit{Lattice congruences of the weak order.}
Order \textbf{21} (2004) no.~4, 315--344. 
 
\bibitem{con_app}
N.~Reading,
\textit{Lattice congruences, fans and Hopf algebras.}
J. Combin. Theory Ser. A {\bf 110} (2005) no.~2, 237--273.

\bibitem{cambrian}
N.~Reading, 
\textit{Cambrian Lattices.}
Adv. Math. \textbf{205} (2006), no.~2, 313--353.

\bibitem{sortable}
N.~Reading,
\textit{Clusters, Coxeter-sortable elements and noncrossing partitions.}
Trans.\ Amer.\ Math.\ Soc.\ \textbf{359} (2007), no.~12, 5931--5958.

\bibitem{sort_camb}
N.~Reading,
\textit{Sortable elements and Cambrian lattices.}
Algebra Universalis \textbf{56} (2007), no.~3-4, 411--437. 

\bibitem{rotate}
N.~Reading,
\textit{Chains in the noncrossing partition lattice.}
SIAM J.\ Discrete Math. \textbf{22} (2008), no.~3, 875-886.

\bibitem{FPSAC}
N.~Reading,
\textit{Noncrossing partitions and the shard intersection order} (Extended abstract).
Formal Power Series and Algebraic Combinatorics (FPSAC), 2009.
\emph{Discrete Math. Theor. Comput. Sci.}, to appear.

\bibitem{camb_fan}
N.~Reading and D.~Speyer,
\emph{Cambrian Fans} (\texttt{math.CO/0606201}).
J. Eur. Math. Soc. (JEMS) \textbf{11} no.~2, 407--447.

\bibitem{typefree}
N. Reading and D. Speyer,
\textit{Sortable elements in infinite Coxeter groups}  (\texttt{arXiv:0803.2722}).
Trans.\ Amer.\ Math.\ Soc., to appear.

\bibitem{ReinerEul}
V.~Reiner,
\textit{The distribution of descents and length in a Coxeter group.}
Electron.\ J.\ Combin.\ \textbf{2} (1995), Research Paper 25, 20 pp. (electronic). 

\bibitem{Rei}
V. Reiner,
\textit{Non-crossing partitions for classical reflection groups},
Discrete Math. \textbf{177} (1997), no.~1-3, 195--222.

\bibitem{Sloane}
N.~J.~A.~Sloane, \textit{The On-Line Encyclopedia of Integer Sequences.}
Published electronically at \texttt{www.research.att.com/\~{}njas/sequences}.

\bibitem{Ziegler}
G.~Ziegler,
\textit{Lectures on polytopes.}
Graduate Texts in Mathematics \textbf{152}. 
Springer-Verlag, New York, 1995.

\end{thebibliography}
\end{document}